\documentclass[11pt]{article}
\usepackage{latexsym}

\newcommand{\eproof}{\mbox{\ }\hfill $\Box$ \par \vskip 10pt}

\newtheorem{Theorem}{Theorem}[section]
\newtheorem{lemma}[Theorem]{Lemma}
\newtheorem{prop}[Theorem]{Proposition}

\newtheorem{corol}[Theorem]{Corollary}

\begin{document}

\title{Transmission eigenvalues for strictly concave domains}

\author{{\sc Georgi Vodev}}

\date{}

\maketitle

\noindent
{\bf Abstract.} We show that for strictly concave domains there are no interior transmission eigenvalues in a region of the form
$\left\{\lambda\in{\bf C}:{\rm Re}\,\lambda\ge 0,\,\,|{\rm Im}\,\lambda|\ge C_\varepsilon\left({\rm Re}\,\lambda+1\right)^{\frac{1}{2}+\varepsilon}\right\}$, 
$C_\varepsilon>0$, for every $0<\varepsilon\ll 1$. As a consequence, we obtain Weyl asymptotics for the number of the transmission eigenvalues with an almost optimal remainder term.

\setcounter{section}{0}
\section{Introduction and statement of results}

Let $\Omega\subset{\bf R}^d$, $d\ge 2$, be a bounded, connected domain with a $C^\infty$ smooth boundary $\Gamma=\partial\Omega$. 
A complex number $\lambda\in {\bf C}$, $\lambda\neq 0$, will be said to be a transmission eigenvalue if the following problem has a non-trivial solution:
$$\left\{
\begin{array}{lll}
\left(\nabla c_1(x)\nabla+\lambda n_1(x)\right)u_1=0 &\mbox{in} &\Omega,\\
\left(\nabla c_2(x)\nabla+\lambda n_2(x)\right)u_2=0 &\mbox{in} &\Omega,\\
u_1=u_2,\,\,\, c_1\partial_\nu u_1=c_2\partial_\nu u_2& \mbox{on}& \Gamma,
\end{array}
\right.
\eqno{(1.1)}
$$
where $\nu$ denotes the exterior Euclidean unit normal to $\Gamma$, $c_j,n_j\in C^\infty(\overline\Omega)$, $j=1,2$ are strictly positive real-valued functions.
Let $f\in C^\infty({\bf R}^d)$ be such that $f<0$ in $\Omega$, $f>0$ in ${\bf R}^d\setminus\Omega$, $df\neq 0$ on $\Gamma$.
Given an Hamiltonian $g\in C^\infty(T^*\Omega)$ of the form
$$g(x,\xi)=\sum_{i,j=1}^d g_{ij}(x)\xi_i\xi_j\ge C|\xi|^2,\quad C>0,$$
the boundary $\Gamma$ will be said to be $g-$ strictly concave (viewed from the interior) iff for any $(x,\xi)$ satisfying
$$f(x)=0,\quad g(x,\xi)=1,\quad \{g,f\}(x,\xi)=0,$$
we have
$$\{g,\{g,f\}\}(x,\xi)>0,$$
where $\{\cdot,\cdot\}$ denotes the Poisson brackets. Set $g_j(x,\xi)=\frac{c_j(x)}{n_j(x)}|\xi|^2$. 
 Our main result is the following

\begin{Theorem} Let $\Gamma$ be $g_j-$ strictly concave, $j=1,2$, and assume either the condition
$$c_1(x)\equiv c_2(x),\quad \partial_\nu c_1(x)\equiv \partial_\nu c_2(x),\quad n_1(x)\neq n_2(x)\quad\mbox{on}\quad\Gamma,\eqno{(1.2)}$$
or the condition 
$$(c_1(x)-c_2(x))(c_1(x)n_1(x)-c_2(x)n_2(x))<0\quad\mbox{on}\quad\Gamma.\eqno{(1.3)}$$
Then, for every $0<\varepsilon\ll 1$ there exists a constant $C_\varepsilon>0$ such that there are no transmission eigenvalues in the region 
$$\left\{\lambda\in{\bf C}:{\rm Re}\,\lambda\ge 0,\,\,|{\rm Im}\,\lambda|\ge C_\varepsilon\left({\rm Re}\,\lambda+1\right)^{\frac{1}{2}+\varepsilon
}\right\}.$$
\end{Theorem}
 
  \noindent
{\bf Remark 1.} It has been proved in \cite{kn:V} that, under the conditions (1.2) and (1.3), for arbitrary domains there are no transmission eigenvalues in
$$\left\{\lambda\in{\bf C}:{\rm Re}\,\lambda\ge 0,\,\,|{\rm Im}\,\lambda|\ge C_\varepsilon\left({\rm Re}\,\lambda+1\right)^{\frac{3}{4}+\varepsilon
}\right\}$$
for every $0<\varepsilon\ll 1$.

  \noindent
{\bf Remark 2.} The assumption that $\Gamma$ is strictly concave does not improve the eigenvalue-free regions in
${\rm Re}\,\lambda<0$. Note that it is proved in \cite{kn:V} that for arbitrary domains there are no transmission eigenvalues
in ${\rm Re}\,\lambda\le -C$ for some constant $C>0$ under the assumption (1.2), and in 
  $$\left\{\lambda\in{\bf C}:{\rm Re}\,\lambda\le 0,\,\,|{\rm Im}\,\lambda|\ge C_N\left(|{\rm Re}\,\lambda|+1\right)^{-N}\right\}$$ 
for any $N>1$ under the assumption (1.3).

 \noindent
{\bf Remark 3.} When the function in the left-hand side of (1.3) is strictly positive, large eigenvalue-free regions have been proved
in \cite{kn:V} for arbitrary domains, which however are worse than the eigenvalue-free regions in the cases considered in the
present paper. It seems that in this case no improvement is possible even if the domain is supposed strictly concave.

 \noindent
{\bf Remark 4.} It has been proved recently in \cite{kn:PV} that the total counting function
$N(r)=\#\{\lambda-{\rm trans.\, eig.}:\,|\lambda|\le r^2\}$, $r>1$, satisfies the asymptotics
$$N(r)=(\tau_1+\tau_2)r^d+O_\varepsilon(r^{d-\kappa+\varepsilon}),\quad\forall\,0<\varepsilon\ll 1,$$
where $0<\kappa\le 1$ is such that there are no transmission eigenvalues in the region
$$\left\{\lambda\in{\bf C}:\,|{\rm Im}\,\lambda|\ge C\left(|{\rm Re}\,\lambda|+1\right)^{1-\frac{\kappa}{2}}\right\},\quad C>0,$$
and $$\tau_j=\frac{\omega_d}{(2\pi)^d}\int_\Omega\left(\frac{n_j(x)}{c_j(x)}\right)^{d/2}dx,$$
$\omega_d$ being the volume of the unit ball in ${\bf R}^d$.

Theorem 1.1 and Remark 4 imply the following

\begin{corol}  Under the conditions of Theorem 1.1, the counting function of the transmission eigenvalues satisfies the asymptotics
$$N(r)=(\tau_1+\tau_2)r^d+O_\varepsilon(r^{d-1+\varepsilon}),\quad\forall\,0<\varepsilon\ll 1.\eqno{(1.4)}$$
\end{corol}

To prove Theorem 1.1 we follow the same strategy as in \cite{kn:V}. We first reduce our problem to a semi-classical one by
putting $h=({\rm Re}\,\lambda)^{-1/2}$, $z=h^2\lambda=1+ih^2{\rm Im}\,\lambda$. Thus we have to show that the operator $T(h,z)=c_1N_1(h,z)-
c_2N_2(h,z)$ is invertible for $|{\rm Im}\,z|\ge h^{1-\varepsilon}$, $0<h\ll 1$, $\forall\, 0<\varepsilon\ll 1$ (see Theorem 7.1),
where $N_j$ is the Dirichlet-to-Neumann (DN) map associated to the operator $h^2\nabla c_j\nabla+zn_j$ (see Section 2 for the
precise definition and the main properties). It is shown in \cite{kn:V} that the operator $T(h,z)$ is invertible in the region
$|{\rm Im}\,z|\ge h^{1/2-\varepsilon}$ for an arbitrary domain $\Omega$. In the present paper we show that this region can be
extended to $|{\rm Im}\,z|\ge h^{1-\varepsilon}$ if $\Gamma$ is strictly concave with respect to both $g_1$ and $g_2$.
To do so, we have to study more carefully the DN map $N_j$ near the glancing manifold $\Sigma_j=\{(x,\xi)\in T^*\Gamma:r_0(x,\xi)=m_j(x)\}$,
where $m_j$ denotes the restriction on $\Gamma$ of the function $n_j/c_j$, while $r_0>0$ is the principal symbol of the Laplace-Beltrami
operator on $\Gamma$ with Riemannian metric induced by the Euclidean metric in ${\bf R}^d$. We show that
$N_j(h,z)={\cal O}(h^{\varepsilon/4}):L^2(\Gamma)\to L^2(\Gamma)$ in an ${\cal O}(h^\varepsilon)$ neighbourhood of $\Sigma_j$ as long as
$h^{1-\varepsilon}\le |{\rm Im}\,z|\le h^{\varepsilon}$ (see Theorem 2.2). With this property in hands, the invertibility of $T$ near $\Sigma_j$ is 
almost immediate since the conditions (1.2) and (1.3) guarantee that $N_{3-j}$ is elliptic on $\Sigma_j$, $j=1,2$. The
invertibility of $T$ outside an ${\cal O}(h^\varepsilon)$ neighbourhood of $\Sigma_1\cup\Sigma_2$ for $|{\rm Im}\,z|\ge h^{1-\varepsilon}$ is much easier and can be done
in precisely the same way as in \cite{kn:V} for an arbitrary domain. Indeed, the conditions (1.2) and (1.3) imply that in this region
$T(h,z)$ is an elliptic $h-\Psi$DO, and hence easy to invert. 

Thus the main (and the most difficult) point in our proof is the estimate (2.7) of Theorem 2.2 concerning the behavior of the DN map near
the glancing manifold. Therefore the present paper is almost entirely devoted to the proof of Theorem 2.2. To do so, we make use
of the global symplectic normal form proved in \cite{kn:PoV} in order to transform our boundary-value problem in an
${\cal O}(h^\varepsilon)$ neighbourhood of the glancing manifold to a much simpler one in which we have complete separation
of the normal and tangential variables (see the model equation in Section 5). The advantage is that we can build a relatively simple parametrix in terms of the Airy function
and its derivatives (see Section 5). Note that our parametrix is much simpler than the parametrix of Melrose-Taylor \cite{kn:MT}
and therefore easier to work with. In particular, it is easier to control it as $|{\rm Im}\,z|\to 0$. Using the properties of the
Airy function (see Section 3) we show in Section 5 that our parametrix is valid in an ${\cal O}(h^{1+\varepsilon}/|{\rm Im}\,z|)$ neighbourhood of the glancing manifold as long as $h^{1-2\varepsilon}\le |{\rm Im}\,z|\le h^{\varepsilon}$. To cover the entire ${\cal O}(h^\varepsilon)$ neighbourhood of the glancing manifold we have to build another parametrix in Section 6 following the parametrix construction in \cite{kn:V} and showing that it
can be improved in the case of our model equation. When $|{\rm Im}\,z|\sim h^{2/3}$ a different parametrix, without using the Airy function, is
constructed by Sj\"ostrand (see Section 11 of \cite{kn:Sj}). In this case, it provides another proof of the estimate (2.7). 
Note finally that in Section 3 we prove some properties of the Airy function which
play a crucial role in the parametrix construction in Section 5. They are more or less well-known and most of them can be found in \cite{kn:O}
and in the appendix of \cite{kn:MT}. 

\section{The Dirichlet-to-Neumann map}

Let $(X,g)$ be a compact Riemannian manifold of dimension $n={\rm dim}\, X\ge 2$ with a non-empty smooth boundary $\partial X$.
Then $(\partial X,\widetilde g)$ is a Riemannian manifold without boundary of dimension $n-1$, where $\widetilde g$ is the Riemannian
metric on $\partial X$ induced by the metric $g$. Denote by $\Delta_X$ and $\Delta_{\partial X}$ the (negative) Laplace-Beltrami operators on
$(X,g)$ and $(\partial X,\widetilde g)$, respectively. The boundary $\partial X$ is said to be strictly concave if the second fundamental form of $\partial X$ is strictly positive. In the case when $X\subset{\bf R}^n$ this definition coincides with that one given in the previous section. Given a function $f\in H^1(\partial X)$, let $u$ solve the equation
$$\left\{
\begin{array}{lll}
 \left(h^2\Delta_X+1+i\mu\right)u=0&\mbox{in}& X,\\
 u=f&\mbox{on}&\partial X,
\end{array}
\right.
\eqno{(2.1)}
$$
where $0<h\ll 1$ is a semi-classical parameter and $\mu\in{\bf R}$, $0<|\mu|\le 1$.
Then the semi-classical Dirichlet-to-Neumann (DN) map 
$$N(h,\mu):H^1(\partial X)\to L^2(\partial X)$$
is defined by
$$N(h,\mu)f:={\cal D}_\nu u|_{\partial X},$$
 where ${\cal D}_\nu=-ih\partial_\nu$, $\nu$ being the unit normal to $\partial X$. It is well-known that for arbitrary
 manifolds one has the bound
 $$\left\|N(h,\mu)\right\|_{H_h^1(\partial X)\to L^2(\partial X)}\le \frac{C}{|\mu|}\eqno{(2.2)}$$
 with a constant $C>0$ independent of $h$ and $\mu$, where $H_h^1(\partial X)$ denotes the Sobolev space 
 $H^1(\partial X)$ equipped with the semi-classical norm $\|(1-h^2\Delta_{\partial X})^{1/2}f\|_{L^2(\partial X)}$.
 It has been proved recently that better bounds are possible if $\mu$ is not too close to zero. Indeed, it follows from
 Theorem 3.2 of \cite{kn:V}, still for arbitrary manifolds, that for every $\varepsilon>0$ there is a constant 
 $0<h_0(\varepsilon)\ll 1$ such that for all $0<h\le h_0$, $|\mu|\ge h^{\frac{1}{2}-\varepsilon}$ we have the bound
 $$\left\|N(h,\mu)\right\|_{H_h^1(\partial X)\to L^2(\partial X)}\le C\eqno{(2.3)}$$
 with a constant $C>0$ independent of $h$, $\mu$ and $\varepsilon$.
 Note that (2.3) does not follow from (2.2).  In \cite{kn:V} semi-classical parametrices of the operator $N(h,\mu)$ are constructed in the hyperbolic zone ${\cal H}=\{(x',\xi')\in
 T^*\partial X:r_0(x',\xi')<1\}$, in the glancing zone ${\cal G}=\{(x',\xi')\in
 T^*\partial X:r_0(x',\xi')=1\}$ and in the elliptic zone ${\cal E}=\{(x',\xi')\in
 T^*\partial X:r_0(x',\xi')>1\}$. Hereafter, $r_0(x',\xi')$ denotes the principal symbol of the operator $-\Delta_{\partial X}$
 written in the coordinates $(x',\xi')$. To be more precise, introduce the set ${\cal S}^{k}_\delta$, $k\in{\bf R}$, $0\le\delta<\frac{1}{2}$,
  of all functions $a\in C^\infty(T^*\partial X)$ satisfying
$$\left|\partial_{x'}^\alpha\partial_{\xi'}^\beta a(x',\xi')\right|\le C_{\alpha,\beta}h^{-\delta(|\alpha|+|\beta|)}
\langle\xi'\rangle^{k-|\beta|}$$
for all multi-indices $\alpha,\beta$ with constants $C_{\alpha,\beta}>0$ independent of $h$. We will denote by ${\rm OP}{\cal S}^{k}_\delta$
the set of the $h$-pseudo-differential operators ($h$-$\Psi$DOs) with symbols in ${\cal S}^{k}_\delta$ defined as follows
 $$\left({\rm Op}_h(a)f\right)(x')=
 \left(\frac{1}{2\pi h}\right)^{n-1}\int_{T^*\partial X}e^{-\frac{i}{h}\langle x'-y',\xi'\rangle}a(x',\xi')f(y')dy'd\xi'.$$
 Let $\chi^-,\chi^0,\chi^+\in C^\infty(T^*\partial X)$ be 
independent of $h$ and such that $\chi^-+\chi^0+\chi^+\equiv 1$, supp$\,\chi^-\subset{\cal H}$, 
supp$\,\chi^+\subset{\cal E}$, $\chi^0$ is supported in a small $h$-independent neighbourhood of ${\cal G}$, 
 $\chi^0=1$ in a smaller $h$-independent neighbourhood of ${\cal G}$. Set $\rho(x',\xi',\mu)=\sqrt{-r_0(x',\xi')+1+i\mu}$ with 
 ${\rm Im}\,\rho>0$. It was shown in \cite{kn:V} that, mod ${\cal O}(h^\infty)$, the operator $N(h,\mu){\rm Op}_h(\chi^-)$ belongs
 to ${\rm OP}{\cal S}^{0}_0$ for $|\mu|\ge h^{1-\varepsilon}$, $0<\varepsilon\ll 1$, with a principal symbol $\rho\chi^-$, 
 the operator $N(h,\mu){\rm Op}_h(\chi^0)$ belongs
 to ${\rm OP}{\cal S}^0_{1/2-\varepsilon}$ for $|\mu|\ge h^{1/2-\varepsilon}$ with a principal symbol $\rho\chi^0$, and 
 $N(h,\mu){\rm Op}_h(\chi^+)$ belongs
 to ${\rm OP}{\cal S}^1_0$ with a principal symbol $\rho\chi^+$. Summing up, we conclude that, mod ${\cal O}(h^\infty)$, the operator $N(h,\mu)$ belongs
 to ${\rm OP}{\cal S}^1_{1/2-\varepsilon}$ for $|\mu|\ge h^{1/2-\varepsilon}$ with a principal symbol $\rho$.
 Therefore, in this case the bound (2.3) is a consequence of well-known properties of the $h$-$\Psi$DOs.
 In fact, a more detailed anaysis of the operator $N(h,\mu)$ can be carried out allowing the functions $\chi^+$, $\chi^-$ and $\chi^0$
 to depend on $h$. More generally, it follows from the analysis in \cite{kn:V} that given any function $\chi\in C_0^\infty(T^*\partial X)$,
 for arbitrary $\partial X$, one can construct a parametrix for the operator  $N(h,\mu){\rm Op}_h(\chi)$ as long as
 $$\min_{{\rm supp}\,\chi}|\rho|^2\ge \frac{h^{1-\varepsilon}}{|\mu|}\quad\mbox{for some}\quad\varepsilon>0.$$
 It is easy to see that 
 given a parameter $0<\delta\ll 1$, there are
 functions $\chi_\delta^-,\chi_\delta^0,\chi_\delta^+\in {\cal S}^{0}_\delta$ such that 
 $\chi_\delta^-+\chi_\delta^0+\chi_\delta^+\equiv 1$, supp$\,\chi_\delta^-\subset\{r_0-1\le -h^\delta\}$, 
supp$\,\chi_\delta^+\subset\{r_0-1\ge h^\delta\}$, supp$\,\chi_\delta^0\subset\{|r_0-1|\le 2h^\delta\}$, 
 $\chi_\delta^0=1$ on $\{|r_0-1|\le h^\delta\}$. As in \cite{kn:V} one can prove the following
 
 \begin{Theorem} For every $0<\varepsilon\ll 1$ there is $h_0(\varepsilon)>0$ such that for $0<h\le h_0$,
 $h^{1-\varepsilon}\le|\mu|\le h^{\varepsilon}$, we have the bound
 $$\left\|N(h,\mu){\rm Op}_h(\chi_{\varepsilon/2}^-)-{\rm Op}_h(\rho\chi_{\varepsilon/2}^-)\right\|_{L^2(\partial X)\to L^2(\partial X)}\le Ch^{1/2}.\eqno{(2.4)}$$
 For $|\mu|\le h^\varepsilon$ we also have the bound
 $$\left\|N(h,\mu){\rm Op}_h(\chi_{\varepsilon/2}^+)-{\rm Op}_h(\rho\chi_{\varepsilon/2}^+)\right\|_{L^2(\partial X)\to L^2(\partial X)}\le Ch^{1/2}.\eqno{(2.5)}$$
 For 
 $h^{1/2-\varepsilon}\le|\mu|\le h^{\varepsilon}$, we have the bound
 $$\left\|N(h,\mu){\rm Op}_h(\chi_{\varepsilon/2}^0)\right\|_{L^2(\partial X)\to L^2(\partial X)}\le Ch^{\varepsilon/4}.\eqno{(2.6)}$$
 \end{Theorem}
 
 When $\partial X$ is strictly concave, Sj\"ostrand showed (see Section 11 of \cite{kn:Sj})
 that (2.3) still holds for $C_1h^{2/3}\le|\mu|\le C_2h^{2/3}$, $C_2>C_1>0$ being arbitrary, independent of $h$ and $\mu$. We will show in the present paper that for strictly concave $\partial X$ the bound (2.3) holds true for $ h^{1-\varepsilon}\le|\mu|\le h^{\varepsilon}$, $\forall \, 0<\varepsilon\ll 1$. 
 To this end, we need to improve only the bound (2.6). We have the following
 
  \begin{Theorem} If $\partial X$ is strictly concave, for every $0<\varepsilon\ll 1$ there is $h_0(\varepsilon)>0$ such that for $0<h\le h_0$,
 $h^{1-\varepsilon}\le|\mu|\le h^{\varepsilon}$, we have the bound
 $$\left\|N(h,\mu){\rm Op}_h(\chi_{\varepsilon/2}^0)\right\|_{L^2(\partial X)\to L^2(\partial X)}\le Ch^{\varepsilon/4}.\eqno{(2.7)}$$
 \end{Theorem}
  
{\it Proof.} We will make use of the symplectic normal form obtained in \cite{kn:PoV} to reduce
our problem to a simpler one for which it is easier to construct a parametrix. This model problem will be studied
in the next sections. Let $y=(y_1,y')\in X_\delta:=(-\delta,\delta)\times\partial X$, $0<\delta\ll 1$, be the normal geodesic coordinates with respect to the Riemannian metric $g$. Here we identify the points in $(0,\delta)\times\partial X$ with
$\{x\in X:{\rm dist}(x,\partial X)<\delta\}$. Then in these coordinates we can write
  $$-h^2\Delta_X={\cal D}_{y_1}^2+q(y_1,y',{\cal D}_{y'})+\mbox{lower order terms},$$
  where ${\cal D}_{y_1}=-ih\partial_{y_1}$, ${\cal D}_{y'}=-ih\partial_{y'}$, $q(y_1,y',\eta')=\sum_{|\alpha|=2}q_\alpha(y_1,y')\eta'^\alpha$.
  Moreover
  $q_0(y',\eta'):=q(0,y',\eta')$
  is the principal symbol of $-\Delta_{\partial X}$ written in the coordinates $(y',\eta')$, while 
  $$q_1(y',\eta'):=\frac{\partial q}{\partial y_1}(0,y',\eta')>0$$
  is the second fundamental form of $\partial X$ supposed to be strictly positive (which is nothing else but the definition of
  $g-$ strictly concavity). Then the principal symbol $p$ of the operator $P(h,\mu)=-h^2\Delta_X-1-i\mu$ can be written in the coordinates
  $(y,\eta)\in T^*X_\delta$ as follows
  $$p(y,\eta)=\eta_1^2+q(y_1,y',\eta')-1-i\mu=\eta_1^2+q_0(y',\eta')+y_1q_1(y',\eta')-1-i\mu+{\cal O}(y_1^2q_0).$$
  Denote by ${\cal R}$ the set of all functions $a\in C^\infty(T^*X_\delta)$ satisfying (with all derivatives)
  $$a={\cal O}(x_1^\infty)+
  {\cal O}(\xi_1^\infty)+{\cal O}((1-q_0)^\infty)$$
  in a neighbourhood of ${\cal K}=\{x_1=\xi_1=1-q_0=0\}$. We will also denote by ${\rm OP}{\cal R}$ the $h-\Psi$DOs on $X_\delta$
  with symbols of the form $\sum_{j=0}^\infty h^ja_j$, where $a_j\in{\cal R}$ do not depend on $h$. Let
  $\phi\in C^\infty({\bf R})$, $\phi(\sigma)=1$ for $|\sigma|\le 1/2$, $\phi(\sigma)=0$ for $|\sigma|\ge 1$. Given
  any $0<\varepsilon\ll 1$, denote by ${\cal A}_\varepsilon$ the $h-\Psi$DO on $X_\delta$ with symbol
  $\phi(x_1/h^\varepsilon)\phi((1-q_0)/h^\varepsilon)$. Clearly, if $R\in {\rm OP}{\cal R}$, we have
  $$R{\cal A}_\varepsilon,\,\,{\cal A}_\varepsilon R = {\cal O}(h^\infty): L^2(X_\delta)\to L^2(X_\delta).$$
  It is shown in \cite{kn:PoV} (see Theorem 3.1) that there exists an exact symplectic map $\chi:T^*X_\delta\to T^*X_\delta$
  such that $\chi(x,\xi)=(y(x,\xi),\eta(x,\xi))$ satisfies
  $$y_1=x_1q_1(x',\xi')^{-1/3}+{\cal O}(x_1^2)+{\cal O}(x_1(1-q_0)),$$ $$\eta_1=\xi_1q_1(x',\xi')^{1/3}+{\cal O}(x_1)+{\cal O}(\xi_1(1-q_0)),$$
  $$(y',\eta')=(x',\xi')+{\cal O}(x_1),$$
  $$(p\circ\chi)(x,\xi)=\left(q_1(x',\xi')^{2/3}+{\cal O}(x_1)\right)(\xi_1^2+x_1-\zeta(x',\xi'))\quad({\mbox mod}\,\,{\cal R})$$
  in a neighbourhood of ${\cal K}$, where 
  $$\zeta(x',\xi')=\left(q_1(x',\xi')^{-2/3}+{\cal O}(1-q_0)\right)(1+i\mu-q_0(x',\xi')).$$
  Thus, if ${\cal U}\subset T^*X_\delta$ is a small neighbourhood of ${\cal K}$, then $\chi$ sends ${\cal U}$ into itself. 
  Using $h-$ Fourier integral operators on $X_\delta$ ($h-$ FIOs) associated to the canonical relation
  $$\Lambda = \{(y,\eta,x,\xi)\in T^*X_\delta\times T^*X_\delta: (y,\eta)=\chi(x,\xi),\,\,(x,\xi)\in{\cal U}\}$$
  one can transform the operator $P$ into a simpler one, $P'_0$, which can be written in the coordinates $(x,\xi)$ as follows
  $$P'_0={\cal D}_{x_1}^2+x_1-L_1(x',{\cal D}_{x'};h)-i\mu L_2(x',{\cal D}_{x'};h)$$
  where $L_j(x',\xi';h)=\sum_{k=0}^\infty h^kL_j^{(k)}(x',\xi')$, $j=1,2,$ with
  $$L_1^{(0)}(x',\xi')=\left(q_1(x',\xi')^{-2/3}+{\cal O}(1-q_0)\right)(1-q_0(x',\xi')),$$ $$ 
  L_2^{(0)}(x',\xi')=q_1(x',\xi')^{-2/3}+{\cal O}(1-q_0).$$
  More precisely, there exist zero-order elliptic (in ${\cal U}$) $h-\Psi$DOs on $X_\delta$, $A$, $A'$, and a
  zero-order elliptic $h-$ FIO on $X_\delta$, $U$, associated to $\Lambda$, such that if we set $T=UA$, $T'=UA'$,
  we have the relations (see Theorem 4.2 of \cite{kn:PoV}):
  $$PT=T'P'_0+T'R_0,\eqno{(2.8)}$$
  $$\iota^*T={\cal Q}_1\iota^*+h{\cal Q}_2\iota^*{\cal D}_{x_1}+\iota^*VP'_0+\iota^*R,\eqno{(2.9)}$$
  $$\iota^*{\cal D}_{x_1}T=\widetilde{\cal Q}_1\iota^*{\cal D}_{x_1}+h\widetilde{\cal Q}_2\iota^*+\iota^*\widetilde VP'_0+\iota^*
  \widetilde R,\eqno{(2.10)}$$
  where $\iota^*$ deontes the restriction on $x_1=0$, ${\cal Q}_j$, $\widetilde{\cal Q}_j$, $j=1,2$, are 
  zero-order $h-\Psi$DOs on $\partial X$, ${\cal Q}_1$ and  $\widetilde{\cal Q}_1$ being elliptic in a neighbourhood
  of $\{q_0=1\}$, 
  $V$ and $\widetilde V$ are zero-order $h-\Psi$DOs on $X_\delta$, and $R_0, R, \widetilde R\in {\rm OP}{\cal R}$.
  One can further simplify the operator $P'_0$ by making a new symplectic change of the tangential variables 
  $(x^\sharp,\xi^\sharp)=\chi_\sharp(x',\xi')\in T^*\partial X$ such that
  $$\xi^\sharp_n=-L_1^{(0)}(x',\xi').$$
  Then, in these coordinates the glancing manifold $\{q_0=1\}$ is defined by $\xi^\sharp_n=0$. 
  Conjugating with a zero-order elliptic (in a neighbourhood of the glancing manifold) $h-$FIO operator on
  $\partial X$ we get (2.8), (2.9) and (2.10) with new operators of the same type (which we will denote in the same way below) and $P'_0$ replaced by
  $$P_0={\cal D}_{x_1}^2+x_1+{\cal D}_{x^\sharp_n}-i\mu Q_0(x^\sharp,{\cal D}_{x^\sharp})+{\cal Q}(x^\sharp,{\cal D}_{x^\sharp}; \mu,h) $$
  where $Q_0(x^\sharp,\xi^\sharp)>0$ in a neighbourhood of $\xi^\sharp_n=0$, and 
  $${\cal Q}=\sum_{k=1}^\infty h^kQ_k(x^\sharp,\xi^\sharp; \mu).$$
  Thus we get the model operator studied in Sections 5 and 6. Indeed, given a function $\widetilde f\in L^2(\partial X)$, it is constructed
   a parametrix $\widetilde u(x_1,x^\sharp)$ supported in $0\le x_1\le h^\varepsilon$ such that 
  $$\widetilde u|_{x_1=0}={\rm Op}_h\left(\phi(\xi^\sharp_n/h^\varepsilon)\right)\widetilde f+{\cal O}(h^\infty)\widetilde f,$$
  $$\left\|P_0\widetilde u\right\|_{H^s((0,\delta)\times\partial X)}\le C_Mh^{M}\left\|\widetilde f\right\|_{L^2(\partial X)}\eqno{(2.11)}$$
  for every $s\ge 0$, where $M\gg 1$ is an arbitrary integer independent of $h$. Hereafter, the Sobolev spaces $H^s$
  will be equipped with the semi-classical norm. Moreover, by Theorem 6.6 the operator defined by
  $$\widetilde N\widetilde f:={\cal D}_{x_1}\widetilde u|_{x_1=0}$$
  satisfies the bound
  $$\left\|\widetilde N\right\|_{L^2(\partial X)\to H^s(\partial X)}\le Ch^{\varepsilon/4}.\eqno{(2.12)}$$
  By (2.8), (2.9) and (2.10) (with $P'_0$ replaced by $P_0$) combined with (2.11) and (2.12) we obtain that the function
  $u=T\widetilde u$ satisfies the bounds
  $$\left\|Pu\right\|_{H^s((0,\delta)\times\partial X)}\le C_Mh^{M}\left\|\widetilde f\right\|_{L^2(\partial X)}\eqno{(2.13)}$$
  $$\left\|u|_{\partial X}-(Q_1+hQ_2\widetilde N)\widetilde f\right\|_{H^s(\partial X)}\le C_Mh^M\left\|\widetilde f\right\|_{L^2(\partial X)}\eqno{(2.14)}$$
  $$\left\|{\cal D}_{x_1}u|_{\partial X}\right\|_{L^2(\partial X)}\le Ch^{\varepsilon/4}\left\|\widetilde f\right\|_{L^2(\partial X)}.\eqno{(2.15)}$$
  Given any function $f\in L^2(\partial X)$, let $v$ solve the equation
  $$\left\{
\begin{array}{lll}
 \left(h^2\Delta_X+1+i\mu\right)v=0&\mbox{in}& X,\\
 v={\rm Op}_h\left(\phi((q_0-1)/h^\varepsilon)\right)f&\mbox{on}&\partial X,
\end{array}
\right.
\eqno{(2.16)}
$$
where the function $\phi$ is as above. Let $\phi_1\in C_0^\infty({\bf R})$ be such that $\phi_1=1$ on supp$\,\phi$. Since $Q_1$ is a zero-order $h-\Psi$DO on $\partial X$, elliptic in a neighbourhood of
$\{q_0=1\}$, thete exists a zero-order $h-\Psi$DO, $Q_1^\flat$, elliptic on $T^*\partial X$, such that 
$(Q_1^\flat)^{-1}={\cal O}(1)$ and $(Q_1^\flat-Q_1){\rm Op}_h\left(\phi_1((q_0-1)/h^\varepsilon)\right)={\cal O}(h^\infty)$
as operators on $H^s(\partial X)$, $s\ge 0$. Set
$$Z=\left[Q_1^\flat+hQ_2\widetilde N,{\rm Op}_h\left(\phi_1((q_0-1)/h^\varepsilon)\right)\right]
={\cal O}(h^{1-\varepsilon}):H^s(\partial X)\to H^s(\partial X).$$
Then, for $h$ small enough the operator $Q_1^\flat+Z$ is invertible on $H^s(\partial X)$
and 
$$\left(Q_1^\flat+Z\right)^{-1}={\cal O}(1):H^s(\partial X)\to H^s(\partial X).$$
Denote by $u$ the parametrix above with
$$\widetilde f={\rm Op}_h\left(\phi_1((q_0-1)/h^\varepsilon)\right)(Q_1^\flat+Z)^{-1}{\rm Op}_h\left(\phi((q_0-1)/h^\varepsilon)\right)f.$$
We have 
$$\left\|\widetilde f\right\|_{L^2(\partial X)}\le {\cal O}(1)\left\|f\right\|_{L^2(\partial X)}$$ 
and
$$(Q_1+hQ_2\widetilde N)\widetilde f=(Q_1^\flat+hQ_2\widetilde N)\widetilde f + {\cal O}(h^\infty)f$$
$$= {\rm Op}_h\left(\phi((q_0-1)/h^\varepsilon)\right)f+Z_1f+{\cal O}(h^\infty)f$$
 where we have put
 $$Z_1={\rm Op}_h\left((1-\phi_1)((q_0-1)/h^\varepsilon)\right)Z(Q_1^\flat+Z)^{-1}{\rm Op}_h\left(\phi((q_0-1)/h^\varepsilon)\right).$$
 We need now the following
 
 \begin{lemma} For small $h$ we have $Z_1={\cal O}(h^\infty):L^2(Y)\to L^2(Y)$. 
 \end{lemma}
 
 {\it Proof.} Given any integer $m\ge 1$ we can write
 $$Z(Q_1^\flat+Z)^{-1}=I-Q_1^\flat(Q_1^\flat+Z)^{-1}$$ 
 $$=I-\sum_{k=0}^mQ_1^\flat(-(Q_1^\flat)^{-1}Z)^k(Q_1^\flat)^{-1}-Q_1^\flat(-(Q_1^\flat)^{-1}Z)^{m+1}(I+(Q_1^\flat)^{-1}Z)^{-1}(Q_1^\flat)^{-1}$$
 where $I$ denotes the identity. 
 Hence, to prove the lemma it suffices to show that
 $${\rm Op}_h\left((1-\phi_1)((q_0-1)/h^\varepsilon)\right)Q_1^\flat(-(Q_1^\flat)^{-1}Z)^k(Q_1^\flat)^{-1}
 {\rm Op}_h\left(\phi((q_0-1)/h^\varepsilon)\right)$$ $$={\cal O}(h^\infty):L^2(Y)\to L^2(Y)\eqno{(2.17)}$$
 for every integer $k\ge 0$, and all functions $\phi$, $\phi_1\in C_0^\infty({\bf R})$ independent of $h$ and such that
 $\phi_1=1$ on supp$\,\phi$. For $k=0$, (2.17) follows from well-known properties of the $h-\Psi$DOs. 
  It is easy also to see that (2.17) with $k=1$ implies
 (2.17) for every $k\ge 1$. On the other hand, to prove (2.17) with $k=1$ it suffices to prove it with $\widetilde N$ in place of
 $Q_1^\flat(-(Q_1^\flat)^{-1}Z)^k(Q_1^\flat)^{-1}$. This property of the operator $\widetilde N$, however, 
 follows from Theorem 6.6.
 \eproof
 
By (2.13), (2.14) and Lemma 2.3, we get
 $$\left\|P(v-u)\right\|_{H^s((0,\delta)\times\partial X)}\le C_Mh^{M}\left\|f\right\|_{L^2(\partial X)}\eqno{(2.18)}$$
 $$\left\|(v-u)|_{\partial X}\right\|_{H^s(\partial X)}\le C_Mh^M\left\|f\right\|_{L^2(\partial X)}\eqno{(2.19)}$$
while (2.15) implies
  $$\left\|{\cal D}_{x_1}u|_{\partial X}\right\|_{L^2(\partial X)}\le Ch^{\varepsilon/4}\left\|f\right\|_{L^2(\partial X)}.\eqno{(2.20)}$$
  Let us see that (2.18),(2.19) and (2.20) imply
  $$\left\|{\cal D}_{x_1}v|_{\partial X}\right\|_{L^2(\partial X)}\le Ch^{\varepsilon/4}\left\|f\right\|_{L^2(\partial X)}.\eqno{(2.21)}$$
  Denote by $G_D$ the self-adjoint Dirichlet realization of the operator $-\Delta_X$ on $L^2(X)$. We have
  $$v-u=E\left((v-u)|_{\partial X}\right)+\left(h^2G_D-i\mu\right)^{-1}P(v-u)$$ 
  $$+\left(h^2G_D-i\mu\right)^{-1}\left(h^2\Delta_X+1+i\mu\right)E\left((v-u)|_{\partial X}\right)$$
  where $E={\cal O}(h^{1/2}):H^{s}(\partial X)\to H^{s+1/2}(X)$, $s\ge 0$, is the extension map, $(Ef)|_{\partial X}=f$,
  $$\|f\|_{H^{s}(\partial X)}\le {\cal O}(h^{-1/2})\|Ef\|_{H^{s+1/2}(X)}.$$
  By (2.18), (2.19), with ${\cal D}_\nu=-ih\partial_\nu$, we have
  $$\left\|{\cal D}_\nu(v-u)\right\|_{L^2(\partial X)}\le Ch^{1/2}\left\|E\left((v-u)|_{\partial X}\right)\right\|_{H^{3/2}(X)}$$
  $$+Ch^{1/2}\left\|\left(h^2G_D-i\mu\right)^{-1}P(v-u)\right\|_{H^{3/2}(X)}$$
  $$+Ch^{1/2}\left\|\left(h^2G_D-i\mu\right)^{-1}\left(h^2\Delta_X+1+i\mu\right)E\left((v-u)|_{\partial X}\right)\right\|_{H^{3/2}(X)}$$
  $$\le C\left(1+|\mu|^{-1}\right)\left\|(v-u)|_{\partial X}\right\|_{H^{1}(\partial X)}
  +Ch^{1/2}|\mu|^{-1}\left\| P(v-u)\right\|_{H^{3/2}(X)}$$ $$\le C_Mh^{ M-1}\|f\|_{L^2(\partial X)}\eqno{(2.22)}$$
  provided $h^{1-\varepsilon}\le|\mu|\le h^{\varepsilon}$, where we have used the coercivity (ellipticity) of the operator $G_D$.
  Taking $M$ big enough we deduce (2.21) from (2.20) and (2.22). Clearly, (2.21) implies (2.7).
  
\section{Some properties of the Airy function}

It is well-known that the Airy function ${\rm Ai}(z)$ is an entire function of order $\frac{3}{2}$ with simple zeros $\{\nu_j\}\subset (-\infty,0)$,
$-\nu_j\sim (3\pi/2)^{2/3}j^{2/3}$,  and satisfying the
equation
$$(\partial_z^2-z){\rm Ai}(z)=0.\eqno{(3.1)}$$
Differentiating (3.1) $k$ times leads to the following equation for the derivatives of the Airy function, 
${\rm Ai}^{(k)}(z)=\frac{d^k{\rm Ai}(z)}{dz^k}$,
$$(\partial_z^2-z){\rm Ai}^{(k)}(z)=k{\rm Ai}^{(k-1)}(z).\eqno{(3.2)}$$
It is also known that the Airy function satisfies the identities
$${\rm Ai}(-z)=e^{i\pi/3}{\rm Ai}_+(z)+e^{-i\pi/3}{\rm Ai}_-(z),\eqno{(3.3)}$$
$${\rm Ai}(-z)^{-1}=c_1^\pm F(-z){\rm Ai}_\pm(z)+c_2^\pm {\rm Ai}'_\pm(z),\eqno{(3.4)}$$
where $c_j^\pm$ are some constants and we have put
$${\rm Ai}_\pm(z)={\rm Ai}(ze^{\pm i\pi/3}),$$
$$F(z)=\frac{{\rm Ai}'(z)}{{\rm Ai}(z)}.$$
The functions ${\rm Ai}$ and ${\rm Ai}_\pm$ satisfy
$${\rm Ai}(z)=\overline{{\rm Ai}(\overline z)},\quad {\rm Ai}_+(z)=\overline{{\rm Ai}_-(\overline z)}.\eqno{(3.5)}$$
In particular, this imples $|{\rm Ai}_+(z)|=|{\rm Ai}_-(z)|$ for real $z$. 
For $|\arg z|<\pi$ we also have the formula
$${\rm Ai}(z)=\exp\left(-\frac{2}{3}z^{3/2}\right)B(z),
\eqno{(3.6)}$$
$$B(z)=\pi^{-1}\int_0^\infty e^{-t^2z^{1/2}}\cos\left(\frac{t^3}{3}\right)dt,$$
where $z^{1/2}$ is taken so that ${\rm Re}\,z^{1/2}>0$, that is,
$$z^{1/2}=|z|^{1/2}\exp\left(i\frac{1}{2}\arg z\right),\quad z^{3/2}=|z|^{3/2}\exp\left(i\frac{3}{2}\arg z\right).$$
Observe that
$${\rm Re}\,z^{1/2}\ge\frac{|{\rm Im}\,z|}{2|z|^{1/2}}.$$
 The function $B$ satisfies the asymptotic expansion
$$B(z)=z^{-1/4}\sum_{\ell=0}^\infty (-1)^\ell b_\ell\xi^{-\ell}\eqno{(3.7)}$$
for $|z|\gg 1$, $|\arg z|\le\pi-\delta$, $0<\delta\ll 1$, where $\xi=\frac{2}{3}z^{3/2}$ and $b_\ell$ are strictly positive real numbers,
$b_0=(2\sqrt{\pi})^{-1}$. In view of (3.6), (3.7) provides an asymptotic expansion for the Airy function ${\rm Ai}(z)$. Moreover (3.7) can be differentiated a finite number of times thus getting an asymptotic expansion
for ${\rm Ai}^{(k)}(z)$. In particular, we get that for $|\arg z|\le\pi-\delta$ the function $F(z)$
has the expansion
$$F(z)=-z^{1/2}\sum_{\ell=0}^\infty \widetilde b_\ell\xi^{-\ell},\quad |z|\gg 1,\eqno{(3.8)}$$
where $\widetilde b_0=1$. Moreover, the function $F^{(k)}(z)=\frac{d^kF(z)}{dz^k}$ has the expansion obtained by differentiating
(3.8) $k$ times. The behaviour of the functions ${\rm Ai}(z)$ and $F(z)$ for $z\in \Lambda_\delta:={\bf C}\setminus \{|\arg z|\le\pi-\delta\}$
is more complicated.

\begin{lemma} For ${\rm Im}\,z\neq 0$ and every integer $k\ge 0$, we have the bound
$$\left|F^{(k)}(z)\right|\le C_k|{\rm Im}\,z|^{-k}\left(|z|^{1/2}+|{\rm Im}\,z|^{-1}
\right).\eqno{(3.9)}$$
\end{lemma}

{\it Proof.} Given any $z\in{\bf C}$ with ${\rm Im}\,z\neq 0$, denote $B(z)=\{w\in {\bf C}:|w-z|\le |{\rm Im}\,z|/2\}$.
Since the function $F$ is analytic on $B(z)$, by the Cauchy theorem we have
$$\left|F^{(k)}(z)\right|\le C_k|{\rm Im}\,z|^{-k}\max_{w\in\partial B(z)}|F(w)|.\eqno{(3.10)}$$
It follows from (3.10) that if (3.9) holds with $k=0$, it holds for all $k$.

Since the function $F(z)$ is analytic at $z=0$, there exists a constant $z_0>0$ such that the bound (3.9) 
holds trivially for $|z|\le z_0$. For $|\arg z|\le\pi-\delta$, $|z|\gg 1$, it follows easily from (3.8).
Therefore, we may suppose that $z_0\le |z|\le z_1$, $z_1>z_0>0$ being constants, or $z\in \Lambda_\delta$, $|z|\gg 1$. 
To deal with the first case we will use the Hadamard factorization theorem. Since the zeros of the Airy function are simple, we
can write
$${\rm Ai}(z)=e^{C_1z+C_2}\prod_{j=1}^\infty\left(1-\frac{z}{\nu_j}\right)e^{\frac{z}{\nu_j}}.$$
Hence we can write the function $F$ in the form
$$F(z)=C_1+\sum_{j=1}^\infty\left((z-\nu_j)^{-1}+\nu_j^{-1}\right).$$
Since $\nu_j$ is real; we have
$$|z-\nu_j|^{-1}\le |{\rm Im}\,z|^{-1},$$
while for $|\nu_j|\ge 2|z|$ we have
$$|z-\nu_j|^{-1}\le 2|\nu_j|^{-1}.$$
Thus we obtain
$$|F(z)|\le |C_1|+\sum_{j=1}^{2|z|}\left(|z-\nu_j|^{-1}+|\nu_j|^{-1}\right)
+|z|\sum_{j=2|z|}^\infty|z-\nu_j|^{-1}|\nu_j|^{-1}$$
$$\le |C_1|+2|z|+2|z||{\rm Im}\,z|^{-1}+2|z|\sum_{j=1}^\infty|\nu_j|^{-2}$$
which gives the desired bound for $|F(z)|$ in this case.

In the second case we will use (3.3). Let $-z\in \Lambda_\delta$, $|z|\gg 1$. Then $|\arg z|\le \delta$ and if $\xi=\frac{2}{3}z^{3/2}$, we have
$${\rm Im}\,\xi={\rm Im}\,z({\rm Re}\,z)^{1/2}(1+{\cal O}(\delta)).$$
Hence
$$|{\rm Im}\,\xi|\ge C_\delta|{\rm Im}\,z||z|^{1/2},\quad C_\delta>0.\eqno{(3.11)}$$
It suffices to consider the case ${\rm Im}\,z>0$ since the case ${\rm Im}\,z<0$ is similar. Then we have ${\rm Im}\,\xi>0$.
In view of (3.7), the functions $B_\pm(z)=z^{1/4}e^{\mp i\pi/12}B(e^{\pm i\pi/3}z)$ satisfy the asymptotics
$$B_\pm(z)=b_0\pm ib_1\xi^{-1}+{\cal O}(\xi^{-2}),\quad -zB'_\pm(z)=\pm \frac{3ib_1}{2}\xi^{-1}+{\cal O}(\xi^{-2}),$$
where $b_0, b_1>0$ are constants. In particular, we have
$$\pm {\rm Im}\,\left(B_\pm(z)\overline{B'_\pm(z)}\right)=\frac{3b_0b_1}{2}|z|^{-5/2}\left(1+{\cal O}(\delta)+{\cal O}(|z|^{-3/2})\right)>0.
\eqno{(3.12)}$$
Let us see that (3.12) implies the inequality
$$|B_+(z)|\ge |B_-(z)|.\eqno{(3.13)}$$
To this end, observe that the first derivative of the function
$$f(\tau)=|B_+({\rm Re}\,z+i\tau)|^2- |B_-({\rm Re}\,z+i\tau))|^2$$
is given by
$$f'(\tau)=2{\rm Im}\,\left(B_+({\rm Re}\,z+i\tau)\overline{B'_+({\rm Re}\,z+i\tau)}\right)-
2{\rm Im}\,\left(B_-({\rm Re}\,z+i\tau)\overline{B'_-({\rm Re}\,z+i\tau)}\right).$$
By (3.12) we get $f'(\tau)>0$ as long as $0\le\tau\le\delta {\rm Re}\,z$ and ${\rm Re}\,z\gg 1$.
On the other hand, in view of (3.5) we have $f(0)=0$. Hence $f(\tau)\ge 0$ for $\tau\ge 0$, which proves (3.13).

By (3.6) and (3.13) we have 
$$\left|\frac{{\rm Ai}_-(z)}{{\rm Ai}_+(z)}\right| = e^{-2{\rm Im}\,\xi}\left|\frac{B_-(z)}{B_+(z)}\right|\le
e^{-2{\rm Im}\,\xi}.\eqno{(3.14)}$$
It is easy to see that the above asymptotics also lead to the bounds
$$\left|\frac{{\rm Ai}'_-(z)}{{\rm Ai}'_+(z)}\right| \le C,\quad \left|\frac{{\rm Ai}'_+(z)}{{\rm Ai}_+(z)}\right| \le C|z|^{1/2}
\eqno{(3.15)}$$
with some constant $C>0$. By (3.11), (3.14) and (3.15),
$$|F(-z)|\le \left|\frac{{\rm Ai}'_+(z)}{{\rm Ai}_+(z)}\right| \left(1+\left|\frac{{\rm Ai}'_-(z)}{{\rm Ai}'_+(z)}\right| \right)
\left(1-\left|\frac{{\rm Ai}_-(z)}{{\rm Ai}_+(z)}\right|\right)^{-1}$$
$$\le\frac{C|z|^{1/2}}{1-e^{-2{\rm Im}\,\xi}}\le\frac{C|z|^{1/2}}{\min\{1,2{\rm Im}\,\xi\}}\le 
C|z|^{1/2}+C|{\rm Im}\,z|^{-1}.$$
\eproof

Given any integer $k\ge 0$, set
$$\Phi_k(z)={\rm Ai}(z)\partial_z^k \left({\rm Ai}(z)^{-1}\right)=\partial_z\Phi_{k-1}(z)-F(z)\Phi_{k-1}(z)\eqno{(3.16)}$$
where $\Phi_{-1}=0$. Clearly, $\Phi_0=1$ and $\Phi_1=-F$. 

\begin{lemma} For ${\rm Im}\,z\neq 0$ and all integers $k\ge 1$, $\ell\ge 0$, we have the bound
$$\left|\partial_z^\ell \Phi_k(z)\right|\le C_{k,\ell}|{\rm Im}\,z|^{-\ell}\left(|z|^{1/2}+|{\rm Im}\,z|^{-1}
\right)^k.\eqno{(3.17)}$$
\end{lemma}

{\it Proof.} Differentiating the identity (3.16) $\ell$ times we get
$$\partial_z^\ell \Phi_k(z)=\partial_z^{\ell+1}\Phi_{k-1}(z)-\sum_{j=0}^\ell c_{\ell, j}F^{(j)}(z)
\partial_z^{\ell-j}\Phi_{k-1}(z).\eqno{(3.18)}$$
It is easy to see by induction in $k$ that (3.17) follows from (3.9).
\eproof

 For $t\ge 0$ and $z\in {\bf C}$, $|\arg\,z|<\pi$, set
 $$\Psi_k(t,z)=\frac{{\rm Ai}^{(k)}(t+z)}{{\rm Ai}(z)},\quad \Psi_k^{(\ell)}(t,z)=\partial_z^\ell \Psi_k(t,z).$$
 
\begin{lemma} For ${\rm Im}\,z\neq 0$ and all integers $k\ge 0$, $\ell\ge 0$, we have the bound
$$\left|\Psi_k^{(\ell)}(0,z)\right|\le C_{k,\ell}\,
|{\rm Im}\,z|^{-\ell}\left(|z|^{1/2}+|{\rm Im}\,z|^{-1}
\right)^{k}.\eqno{(3.19)}$$
For $t> 0$, ${\rm Im}\,z\neq 0$ and all integers $k\ge 0$, $\ell\ge 0$, we have the bound
$$\left|\Psi_k^{(\ell)}(t,z)\right|\le C_{k,\ell}\,
|{\rm Im}\,z|^{-\ell}\left(|z|^{1/2}+|{\rm Im}\,z|^{-1}
\right)^{k+1}\eqno{(3.20)}$$
while for $t\ge |z|$ we have
$$\left|\Psi_k^{(\ell)}(t,z)\right|\le C_{k,\ell}\,|{\rm Im}\,z|^{-\ell}\left(|z|^{1/2}+|{\rm Im}\,z|^{-1}
\right)\left(t^{1/2}+|{\rm Im}\,z|^{-1}
\right)^{k}e^{-t^{1/2}|{\rm Im}\,z|/4}.\eqno{(3.21)}$$
\end{lemma}

{\it Proof.} In view of (3.10) with $\Psi_k$ in place of $F$, it suffices to prove these bounds with $\ell=0$. 
Furthermore, using (3.2) it is easy to see by induction in $k$ that (3.9) implies the estimate
$$\left|{\rm Ai}^{(k)}(z)\right|\le C_k\left(|z|^{1/2}+|{\rm Im}\,z|^{-1}
\right)^{k}\left|{\rm Ai}(z)\right|.\eqno{(3.22)}$$
Hence
$$\left|\Psi_k(t,z)\right|\le C_k\left(t^{1/2}+|z|^{1/2}+|{\rm Im}\,z|^{-1}
\right)^{k}\left|\Psi_0(t,z)\right|.\eqno{(3.23)}$$
In particular, (3.23) implies that (3.19) and (3.21) with $\ell=0$, $k\ge 1$, follows from (3.19) and (3.21) with $\ell=0$, $k=0$.
The same conclusion is still valid concerning the bound (3.20) as long as $t\le 2|z|$. For $t\ge 2|z|$, 
(3.20) follows from (3.21) in view of the inequality
$$t^{k/2}e^{-t^{1/2}|{\rm Im}\,z|/4}\le C_k|{\rm Im}\,z|^{-k}.$$
Therefore, to prove the lemma we have to bound $|\Psi_0|$. Clearly, $\Psi_0(0,z)=1$ which proves (3.19). To bound 
$|\Psi_0(t,z)|$ for $t>0$, let us see that the Airy function satisfies the bounds
$$|{\rm Ai}(z)|\le C\langle z\rangle^{-1/4}e^{-\frac{2}{3}{\rm Re}\,z^{3/2}}, \eqno{(3.24)}$$
$$|{\rm Ai}(z)|^{-1}\le C\langle z\rangle^{-1/4}\left(|z|^{1/2}+|{\rm Im}\,z|^{-1}
\right)e^{\frac{2}{3}{\rm Re}\,z^{3/2}}. \eqno{(3.25)}$$
Indeed, for $|\arg z|\le\pi-\delta$, (3.24) and (3.25) follow from (3.6) and (3.7), while for $z\in\Lambda_\delta$
they follow from (3.3) and (3.4) combined with Lemma 3.1. By (3.24) and (3.25),
$$|\Psi_0(t,z)|\le C\left(|z|^{1/2}+|{\rm Im}\,z|^{-1}
\right)e^{-\varphi(t,z)} \eqno{(3.26)}$$
where
$$\varphi=\frac{2}{3}{\rm Re}\,(z+t)^{3/2}-\frac{2}{3}{\rm Re}\,z^{3/2}=\int_0^t{\rm Re}\,(z+\tau)^{1/2}d\tau$$ $$
\ge \frac{1}{2}\int_0^t\frac{|{\rm Im}\,z|}{|z+\tau|^{1/2}}d\tau\ge\frac{t|{\rm Im}\,z|}{2|z|^{1/2}+2t^{1/2}}. \eqno{(3.27)}$$
Hence $\varphi\ge 0$ for $t\ge 0$, while for $t\ge |z|$ we have $\varphi\ge \frac{1}{4}t^{1/2}|{\rm Im}\,z|$.
Therefore, the desired bounds for $|\Psi_0|$ follow from (3.26).
\eproof

  \section{Some properties of the $h-\Psi$ DOs}

Let $Y$ be an $n-1$ dimensional compact manifold without
boundary or an open neighbourhood in ${\bf R}^{n-1}$. In this section we will recall some useful criteria on a
symbol $a'y,\eta)\in T^*Y$ for the $h-\Psi$ DO, ${\rm Op}_h(a)$, to be bounded on $L^2(Y)$. We will make use of the analysis developed
in Section 7 of \cite{kn:DS} (see also Section 2 of \cite{kn:V}). We first have the following

\begin{prop} Let $a\in T^*Y$ satisfy the bounds
$$\left|\partial_y^\alpha a(y,\eta)\right|\le a_0(h)h^{-|\alpha|/2}\eqno{(4.1)}$$
for $|\alpha|\le n$, where $a_0>0$ is a parameter.
Then there is a constant $C>0$ independent of $h$ such that
$$\left\|{\rm Op}_h(a)\right\|_{L^2(Y)\to L^2(Y)}\le Ca_0(h).\eqno{(4.2)}$$
\end{prop}

This proposition follows for example from Proposition 2.1 of \cite{kn:V}. The next proposition can be derived from the analysis in 
Section 7 of \cite{kn:DS}.

\begin{prop} Let $a,b\in T^*Y$ satisfy the bounds
$$\left|\partial_y^\alpha\partial_\eta^\beta a(y,\eta)\right|\le C_{\alpha,\beta},\eqno{(4.3)}$$
$$\left|\partial_y^\alpha b(y,\eta)\right|\le C_{\alpha}h^{-M_0-\delta|\alpha|}\eqno{(4.4)}$$
where $0\le\delta<1$, for all multi-indices $\alpha$ and $\beta$ with constants $C_{\alpha}, C_{\alpha,\beta}>0$ independent of $h$,
and $M_0>0$ independent of $h$ and $\alpha$. Then for every integer $M\gg M_0$ there is a constant $C_M>0$ independent of $h$ such that
$$\left\|{\rm Op}_h(ab)-{\rm Op}_h\left(\sum_{|\alpha|=0}^M\frac{(-ih)^{|\alpha|}}{|\alpha|!}\partial_\eta^\alpha a\partial_y^\alpha b\right)\right\|_{L^2(Y)\to L^2(Y)}\le C_Mh^{M(1-\delta)/2}.\eqno{(4.5)}$$
\end{prop}

{\it Proof.} In view of formula (7.15) 
of \cite{kn:DS} the operator in the left-hand side of (4.5)
whose norm we would like to bound is an $h$-psdo with symbol $c(x,\xi,x,\xi)$, where the function $c$ is given by
$$c(x,\xi,y,\eta)=e^{ih D_\xi\cdot D_y}a(x,\xi)b(y,\eta)-\sum_{|\alpha|=0}^M\frac{(-ih)^{|\alpha|}}{|\alpha|!}\partial_\eta^\alpha a(x,\xi)\partial_y^\alpha b(y,\eta)$$
where we have put $D=-i\partial$. The inequality (7.17) of \cite{kn:DS} together with (4.3) and (4.4) yield the estimate
$$|c(x,\xi,y,\eta)|\le C_{s,M}h^M\sum_{|\alpha|+|\beta|\le s}\left\|D_\xi^\alpha D_y^\beta (D_\xi\cdot D_y)^Ma(x,\xi)b(y,\eta)\right\|_{L^2}$$
 $$\le C_{s,M}h^{M(1-\delta)-M_0-s\delta}\eqno{(4.6)}$$
 for $s>(n-1)/2$. Similarly, for all multi-indices $\alpha$ and $\beta$, we have
 $$|\partial_x^\alpha\partial_y^\beta c(x,\xi,y,\eta)|\le C_{s,M,\alpha,\beta}h^{M(1-\delta)-M_0-s\delta-|\beta|\delta}.\eqno{(4.7)}$$
 By (4.7) we get
 $$|\partial_x^\alpha c(x,\xi,x,\xi)|\le C_{s,M,\alpha}h^{M(1-\delta)-M_0-s\delta-|\alpha|\delta}.\eqno{(4.8)}$$
 By Proposition 4.1 and (4.8), with some $\ell>0$ depending only on the dimension, we conclude
$$\left\|{\rm Op}_h(c(x,\xi,x,\xi))\right\|_{L^2\to L^2}\le C_Mh^{M(1-\delta)-M_0-\ell\delta}\le C_Mh^{M(1-\delta)/2}\eqno{(4.9)}$$
if $M$ is taken large enough. 
\eproof

\section{Parametrix construction for the model equation}

Let the parameters $h$ and $\mu$ be as in Section 2, $h^{1-2\varepsilon}
\le|\mu|\le h^{\varepsilon}$, $0<\varepsilon\ll 1$. Let also $Y$ be as in Section 4. Consider the operator
$$P_0={\cal D}_t^2+t+{\cal D}_{y_1}+i\mu q(y,{\cal D}_y)+h\widetilde q(y,{\cal D}_y;h,\mu),\quad t>0,$$
where ${\cal D}_t=-ih\partial_t$, ${\cal D}_y=-ih\partial_y$, $y\in Y$, the function $q\in C^\infty(T^*Y)$, $q\in{\cal S}_0^0$, is real-valued and does not depend on $t$, $h$ and $\mu$, satisfying $0<C_1\le q\le C_2$, $C_1$ and $C_2$ being constants,
$\widetilde q\in{\cal S}_0^0$ uniformly in $h$ and $\mu$. Let $\eta=(\eta_1,\eta')$ be the 
dual variables of $y=(y_1,y')$. Let also the function $\phi$ be as in Section 2.
We are going to build a parametrix, $\widetilde u$, for the solution $u$ of the equation
$$\left\{
\begin{array}{lll}
 P_0u=0&\mbox{in}& {\bf R}^+\times Y,\\
 u=f_1&\mbox{on} & Y,
\end{array}
\right.
\eqno{(5.1)}
$$
where $f_1$ is microlocally suppoted in the region ${\cal G}(\varepsilon):=\{(\mu,\eta_1)\in {\bf R}^2:|\mu|+|\eta_1|\le 2h^{\varepsilon}\}$. 
We will first construct a parametrix
in the region
$${\cal G}_1(\varepsilon):=\{(\mu,\eta_1)\in {\bf R}^2:|\mu|\left(|\mu|+|\eta_1|\right)\le h^{1+\varepsilon}\}.\eqno{(5.2)}$$
More precisely, in this section we will construct a parametrix, $\widetilde u_1$, of the solution of the equation (5.1) with 
 $f_1={\rm Op}_h\left(\phi(\eta_1|\mu|/h^{1+\varepsilon})\right)f+{\cal O}(h^\infty)f$, $f\in L^2(Y)$ being arbitrary. 
 The construction in the region ${\cal G}_2(\varepsilon):=\{(\mu,\eta_1)\in {\bf R}^2:h^{1+\varepsilon}/|\mu|\le |\mu|+|\eta_1|\le 2h^{\varepsilon}\}$
 will be carried out in the next section.
 
We will be looking for $\widetilde u_1$ in the form
$$\widetilde u_1=\phi(t/h^{\varepsilon}){\rm Op}_h(A(t))g$$
where $g\in L^2(Y)$ will be determined later on such that $\|g\|_{L^2(Y)}\le {\cal O}(1)\|f\|_{L^2(Y)}$, and 
$$A(t)=\sum_{k=0}^M a_k(y,\eta;h,\mu)\psi_k(t,y,\eta;h,\mu),$$
$$\psi_k=h^{k/3}\Psi_k\left(th^{-2/3}, (\eta_1+i\mu q(y,\eta))h^{-2/3}\right),$$
 $\Psi_k$ being the functions introduced in Section 3, $M$ is an arbitrary integer, $a_0=\phi_1(\eta_1|\mu|/h^{1+\varepsilon})$, $\phi_1\in C_0^\infty({\bf R})$ being such that $\phi_1=1$ on supp$\,\phi$, while $a_k$, $k\ge 1$, do not depend on the variable 
$t$ and will be determined later on. Observe first that we have
$$P_0{\rm Op}_h(A(t))={\rm Op}_h\left(({\cal D}_t^2+t+\eta_1+i\mu q(y,\eta)-ih\partial_{y_1})A(t)\right)$$
$$+i\mu q(y,{\cal D}_y){\rm Op}_h(A(t))-i\mu {\rm Op}_h(qA(t))+h\widetilde q(y,{\cal D}_y){\rm Op}_h(A(t)).\eqno{(5.3)}$$
It is easy to see that (3.2) implies the identity
$$({\cal D}_t^2+t+\eta_1+i\mu q(y,\eta))\Psi_k\left(th^{-2/3}, (\eta_1+i\mu q(y,\eta))h^{-2/3}\right)$$ $$=
-kh^{2/3}\Psi_{k-1}\left(th^{-2/3}, (\eta_1+i\mu q(y,\eta))h^{-2/3}\right)$$
and hence
$$({\cal D}_t^2+t+\eta_1+i\mu q(y,\eta))A(t)=-h\sum_{k=0}^{M-1}(k+1)a_{k+1}\psi_k.\eqno{(5.4)}$$
Using the identity
$$\partial_z\Psi_k(z)=\Psi_{k+1}(t,z)-F(z)\Psi_k(t,z)$$
 we can also write
$$\partial_{y_1}\Psi_k\left(th^{-2/3}, (\eta_1+i\mu q(y,\eta))h^{-2/3}\right)=i\mu h^{-2/3}\partial_{y_1}q
\Psi_{k+1}\left(th^{-2/3}, (\eta_1+i\mu q(y,\eta))h^{-2/3}\right)$$
$$-i\mu h^{-2/3}\partial_{y_1}qF\left(\eta_1+i\mu q(y,\eta))h^{-2/3}\right)\Psi_k\left(th^{-2/3}, (\eta_1+i\mu q(y,\eta))h^{-2/3}\right).$$
Hence
$$\partial_{y_1}A(t)=\sum_{k=0}^M\left(\partial_{y_1}a_k-i\mu h^{-1}\partial_{y_1}qF^\sharp a_k+
i\mu h^{-1}\partial_{y_1}qa_{k-1}\right)\psi_k$$ $$+i
\mu h^{-1}\partial_{y_1}qa_{M}\psi_{M+1}\eqno{(5.5)}$$
where $a_{-1}=0$ and we have put
$$F^\sharp=h^{1/3}F\left(\eta_1+i\mu q(y,\eta))h^{-2/3}\right).$$
Set
$$\rho_1=|\eta_1|^{1/2}+|\mu|^{1/2}+\frac{h}{|\mu|}<1.$$

\begin{lemma} For $t= 0$, all $k\ge 0$ and multi-indices $\alpha$, we have the bound
$$\left|\partial_y^\alpha\psi_k\right|\le C_{k,\alpha}\,\rho_1^k.\eqno{(5.6)}$$
For all $t> 0$, $k\ge 0$ and multi-indices $\alpha$, we have the bound
$$\left|\partial_y^\alpha\psi_k\right|\le C_{k,\alpha}\,h^{-1/3}\rho_1^k.\eqno{(5.7)}$$
Moreover, there exists a constant $C>0$ such that for $C(|\mu|+|\eta_1|)\le t\le 1$ we have the bound
$$\left|\partial_y^\alpha\psi_k\right|\le C_{k,\alpha}h^{-1/3}e^{-t^{1/2}|\mu|/4h}.\eqno{(5.8)}$$
We also have the bound
$$\left|\partial_y^\alpha F^\sharp\right|\le C_{\alpha}\,\rho_1.\eqno{(5.9)}$$
\end{lemma}

{\it Proof.}  It is easy to see by induction that
$$\partial_y^\alpha\Psi_k\left(th^{-2/3}, (\eta_1+i\mu q(y,\eta))h^{-2/3}\right)$$ $$=\sum_{j=0}^{|\alpha|}c_{\alpha,j}(y,\eta)\left(\frac{\mu}{h^{2/3}}\right)^j
\Psi_k^{(j)}\left(th^{-2/3}, (\eta_1+i\mu q(y,\eta))h^{-2/3}\right)\eqno{(5.10)}$$
with some function $c_{\alpha,j}$ independent of $t$, $h$ and $\mu$, $c_{\alpha,0}=0$ for $|\alpha|\ge 1$. Recall that $q\ge C_1>0$. Now (5.6)-(5.8) follow from Lemma 3.3 and (5.10). 
The bound (5.9) follows from (3.9) and (5.10) applied with $F^\sharp$ in place of $\Psi_k$.
\eproof

Set
$$E_1(t)=\frac{i\mu}{h}\sum_{|\alpha|=1}^M\frac{(-ih)^{|\alpha|}}{|\alpha|!}\partial_\eta^\alpha q\partial_y^\alpha A(t),$$
$$E_2(t)=\sum_{|\alpha|=0}^M\frac{(-ih)^{|\alpha|}}{|\alpha|!}\partial_\eta^\alpha \widetilde q\partial_y^\alpha A(t),$$
$${\cal E}_1(t)=i\mu\,q(y,{\cal D}_y){\rm Op}_h(A(t))-i\mu\,{\rm Op}_h(qA(t))-h\,{\rm Op}_h\left(E_1(t)\right),$$
$${\cal E}_2(t)=h\widetilde q(y,{\cal D}_y){\rm Op}_h(A(t))-h\,{\rm Op}_h\left(E_2(t)\right).$$

\begin{lemma} We have the identities
$$E_j(t)=\sum_{k=0}^{2M}\sum_{\ell=0}^k\sum_{|\alpha|=0}^k b_{k,\ell,\alpha}^{(j)}(y,\eta;h,\mu)\partial_y^\alpha a_\ell \psi_k\eqno{(5.11)}$$
where the functions $b_{k,\ell,\alpha}^{(j)}$ do not depend on $a_\nu$, $\psi_\nu$, and satisfy the bounds 
$\partial_y^\beta b_{k,\ell,\alpha}^{(j)}={\cal O}_\beta (1)$ for all multi-indices $\beta$ uniformly in $\mu$ and $h$.
\end{lemma}

{\it Proof.} Using the identity
$$\Psi_k^{(\ell)}(t,z)=\sum_{\nu=0}^\ell \gamma_{\ell,\nu}\partial_z^{\ell-\nu}\left({\rm Ai}(z)^{-1}\right){\rm Ai}^{(k+\nu)}(t+z)$$
 $$=\sum_{\nu=0}^\ell \gamma_{\ell,\nu}\Phi_{\ell-\nu}(z)\Psi_{k+\nu}(t,z)$$
 together with (5.10), we get the identity, 
 $$\partial_y^\alpha\psi_k=\sum_{j=0}^{|\alpha|}\sum_{\nu=0}^j\widetilde c_{\alpha,j,\nu}(y,\eta)\left(\frac{\mu}{h}\right)^j
 \Phi^\sharp_{j-\nu}\psi_{k+\nu}\eqno{(5.12)}$$
 where we have put 
 $$\Phi^\sharp_k=h^{k/3}\Phi_k\left((\eta_1+i\mu q(y,\eta))h^{-2/3}\right).$$
 As in the proof of Lemma 5.1, one can deduce from Lemma 3.2 that $\partial_y^\beta\Phi^\sharp_k={\cal O}_{k,\beta}(1)$.
 Therefore, using (5.12) we can write
 $$h^{|\alpha|}\partial_y^\alpha A(t)=\sum_{k=0}^M\sum_{|\alpha_1|+|\alpha_2|=|\alpha|}\gamma_{\alpha_1,\alpha_2}(h\partial_y)^{\alpha_1}a_k
 (h\partial_y)^{\alpha_2}\psi_k$$
 $$=\sum_{k=0}^{M+|\alpha|}\sum_{\ell=0}^k\sum_{|\alpha_1|=0}^{|\alpha|} e_{k,\ell,\alpha_1}(y,\eta;h,\mu)\partial_y^{\alpha_1} a_\ell \psi_k\eqno{(5.13)}$$
 with functions $e_{k,\ell,\alpha_1}$ independent of $a_k$, $\psi_k$, and satisfying the bounds $\partial_y^\beta e_{k,\ell,\alpha_1}={\cal O}_{\beta}(1)$. Moreover, when $|\alpha|\ge 1$ we have $\widetilde c_{\alpha,j,\nu}=0$ for $j=0$ in (5.12), and hence in this case
 $\partial_y^\beta e_{k,\ell,\alpha_1}={\cal O}_{\beta}(|\mu|)$. Since (5.2) implies $|\mu|^2\le h$, 
 it is easy to see that (5.13) implies (5.11).
 \eproof
 
 We let now the functions $a_k$ satisfy the equations
 $$(k+1)a_{k+1}=-i\partial_{y_1}a_k+\mu h^{-1}\partial_{y_1}qF^\sharp a_k-
\mu h^{-1}\partial_{y_1}qa_{k-1}$$
$$+\sum_{\ell=0}^k\sum_{|\alpha|=0}^k \left(b_{k,\ell,\alpha}^{(1)}+b_{k,\ell,\alpha}^{(2)}\right)\partial_y^\alpha a_\ell .\eqno{(5.14)}$$
Set
$$\rho_2=\frac{|\mu|\rho_1}{h}+\sqrt{\frac{|\mu|}{h}}>1.$$
 
 \begin{lemma} For all integers $k\ge 0$ and all multi-indices $\alpha$, we have the bound
$$\left|\partial_y^\alpha a_k\right|\le C_{k,\alpha}\,\rho_2^k.\eqno{(5.15)}$$
 \end{lemma}
 
 {\it Proof.} In view of Lemmas 5.1 and 5.2, differentiating (5.14) we get
 $$\partial_y^\alpha a_{k+1}=\sum_{|\alpha_1|=0}^{|\alpha|+1}{\cal O}(\rho_2^2)\partial_y^{\alpha_1} a_{k-1}+
 \sum_{|\alpha_2|=0}^{|\alpha|}{\cal O}(\rho_2)\partial_y^{\alpha_2} a_{k}
 +\sum_{\ell=0}^k\sum_{|\beta|=0}^{k+|\alpha|} {\cal O}(1)\partial_y^\beta a_\ell .\eqno{(5.16)}$$
 Since (5.15) is trivially fulfilled for $k=0$, it is easy to see by induction in $k$ that (5.16) implies (5.15) for all $k$.
 \eproof
 
 With this choice of the functions $a_k$ the identity (5.3) becomes
 $$P_0{\rm Op}_h(A(t))={\rm Op}_h\left(B(t)\right)+{\cal E}_1(t)+{\cal E}_2(t)\eqno{(5.17)}$$
 where
 $$B(t)=h(M+1)a_{M+1}\psi_M+\mu\partial_{y_1}qa_{M}\psi_{M+1}$$ $$+\sum_{j=1}^2\sum_{k=M+1}^{2M}\sum_{\ell=0}^k\sum_{|\alpha|=0}^k hb_{k,\ell,\alpha}^{(j)}(y,\eta;h,\mu)\partial_y^\alpha a_\ell \psi_k.$$
 Combining Lemmas 5.1, 5.2 and 5.3 leads to the following

\begin{lemma} For $t= 0$, all $k\ge 0$ and multi-indices $\alpha$, we have the bound
$$\left|\partial_y^\alpha(a_k\psi_k)\right|\le C_{k,\alpha}\,(\rho_1\rho_2)^k.\eqno{(5.18)}$$
For all $t\ge 0$, $k\ge 0$ and multi-indices $\alpha$, we have the bounds
$$\left|\partial_y^\alpha(a_k\psi_k)\right|\le C_{k,\alpha}\,h^{-1/3}(\rho_1\rho_2)^k,\eqno{(5.19)}$$
$$\left|\partial_y^\alpha B(t)\right|\le C_{M,\alpha}(\rho_1\rho_2)^M,\eqno{(5.20)}$$
Moreover, there exists a constant $C>0$ such that for $C(|\mu|+|\eta_1|)\le t\le 1$ we have the bound
$$\left|\partial_y^\alpha(a_k\psi_k)\right|\le C_{k,\alpha}h^{-1/3}\rho_2^ke^{-t^{1/2}|\mu|/4h}.\eqno{(5.21)}$$
\end{lemma}
 
 Observe now that the condition (5.2) implies
 $$\rho_1\rho_2\le C\sqrt{\frac{h}{|\mu|}}+C\left(\frac{|\mu|}{h}(|\mu|+|\eta_1|)\right)^{1/2}$$ $$
 +C\frac{h}{|\mu|}+C\frac{|\mu|}{h}(|\mu|+|\eta_1|)\le{\cal O}(h^{\varepsilon/2}).\eqno{(5.22)}$$
 Using Lemma 5.4 together with (5.22) we will prove the following
 
 \begin{prop} For all $s\ge 0$, we have the bounds
 $$\left\|P_0\widetilde u_1\right\|_{H^s({\bf R}^+\times Y)}\le C_{s,M}h^{M\varepsilon/2}\|g\|_{L^2(Y)},\eqno{(5.23)}$$
  $$\left\|{\rm Op}_h(A(0))g-{\rm Op}_h(a_0)g\right\|_{L^2(Y)}\le Ch^{\varepsilon/2}\|g\|_{L^2(Y)},\eqno{(5.24)}$$
   $$\left\|{\rm Op}_h({\cal D}_tA(0))g\right\|_{L^2(Y)}\le Ch^{\varepsilon}\|g\|_{L^2(Y)}.\eqno{(5.25)}$$
 \end{prop}
 
 {\it Proof}. In view of (5.17) we can write
 $$P_0\widetilde u_1=\phi(t/h^\varepsilon)\left({\rm Op}_h\left(B(t)\right)+{\cal E}_1(t)+{\cal E}_2(t)\right)g
 +\left[{\cal D}_t^2,\phi(t/h^\varepsilon)\right]{\rm Op}_h\left(A(t)\right)g.\eqno{(5.26)}$$
 By (5.19) we have $\partial_y^\alpha{\cal D}_t^\ell A(t)={\cal O}_{\alpha,\ell}\left(h^{-1/3}\right)$, $\forall \alpha,\ell$,
 and hence by Proposition 4.2 we get the bound
 $$\left\|\partial_y^\alpha{\cal D}_t^\ell{\cal E}_j(t)g\right\|_{L^2({\bf R}^+\times Y)}\le C_{M,\alpha,\ell}\,h^M\|g\|_{L^2(Y)}.\eqno{(5.27)}$$
 By (5.20) and (5.22) we have $\partial_y^\alpha{\cal D}_t^\ell B(t)={\cal O}_{\alpha,\ell}\left(h^{M\varepsilon/2}\right)$, $\forall \alpha,\ell$,
 and hence by Proposition 4.1 we get the bound
 $$\left\|\partial_y^\alpha{\cal D}_t^\ell {\rm Op}_h\left(B(t)\right)g\right\|_{L^2({\bf R}^+\times Y)}\le C_{M,\alpha,\ell}\,h^{M\varepsilon/2}\|g\|_{L^2(Y)}.\eqno{(5.28)}$$
 On the other hand, since (5.2) implies $|\mu|+|\eta_1|\le h^{2\varepsilon}$, taking $h$ small enough we can arrange that
 $t\ge C(|\mu|+|\eta_1|)$ as long as $t\in{\rm supp}\left[{\cal D}_t^2,\phi(t/h^\varepsilon)\right]$. Therefore, we can use
 (5.21) to conclude that for $t\sim h^\varepsilon$ we have the bounds $\partial_y^\alpha{\cal D}_t^\ell A(t)={\cal O}_{\alpha,\ell}\left(e^{-ch^{-\varepsilon/2}}\right)$, $\forall \alpha,\ell$, with some constant $c>0$.
 Thus, Proposition 4.1 yields the bound
 $$\left\|\partial_y^\alpha{\cal D}_t^\ell \left[{\cal D}_t^2,\phi(t/h^\varepsilon)\right]{\rm Op}_h\left(A(t)\right)g\right\|_{L^2({\bf R}^+\times Y)}\le C_{\alpha,\ell}\,e^{-ch^{-\varepsilon/2}}\|g\|_{L^2(Y)}.\eqno{(5.29)}$$
 Now (5.23) follows from (5.26)-(5.29) by taking $M$ big enough, depending on $\varepsilon$. Since $\psi_0=1$ for $t=0$, the bound (5.24) follows from (5.18), (5.22) and Proposition 4.1. The proof of (5.25) is similar,
 in view of the identity
 $$h\partial_tA(t)=\sum_{k=0}^M a_k\psi_{k+1}.\eqno{(5.30)}$$
 Indeed, by (5.6), (5.15), (5.22) and (5.30), we have 
 $\partial_y^\alpha{\cal D}_tA(0)={\cal O}_{\alpha}\left(\rho_1\right)$, $\forall \alpha$. Therefore, since $\rho_1={\cal O}(h^\varepsilon)$,
 we get (5.25) by Proposition 4.1.
 \eproof
 
 Set $Z={\rm Op}_h(A(0)-a_0)$. Since the estimate (5.24) holds for every $g\in L^2(Y)$, we have $Z={\cal O}(h^{\varepsilon/2}):L^2(Y)\to L^2(Y)$.
 Hence 
 the operator $I+Z$ is invertible on $L^2(Y)$ for small $h$. Given any $f\in L^2(Y)$, take now
 $$g=(I+Z)^{-1}{\rm Op}_h\left(\phi(\eta_1|\mu|/h^{1+\varepsilon})\right)f.$$
 With this choice of $g$ we have
 $$\widetilde u_1|_{t=0}={\rm Op}_h(A(0))g={\rm Op}_h\left(\phi(\eta_1|\mu|/h^{1+\varepsilon})\right)f+Z_1f$$ 
 where we have put
 $$Z_1=
 {\rm Op}_h\left((1-\phi_1)(\eta_1|\mu|/h^{1+\varepsilon})\right)(I+Z)^{-1}{\rm Op}_h\left(\phi(\eta_1|\mu|/h^{1+\varepsilon})\right).$$
 Thus, to complete the parametrix construction in this case we have to prove the following
 
 \begin{lemma} For small $h$ we have $Z_1={\cal O}(h^\infty):L^2(Y)\to L^2(Y)$. 
 \end{lemma}
 
 {\it Proof.} Given any integer $m\ge 1$ we can write
 $$(I+Z)^{-1}=\sum_{k=0}^m(-Z)^k+(-Z)^{m+1}(I+Z)^{-1}.$$
 Hence, to prove the lemma it suffices to show that
 $${\rm Op}_h\left((1-\phi_1)(\eta_1|\mu|/h^{1+\varepsilon})\right)Z^k{\rm Op}_h\left(\phi(\eta_1|\mu|/h^{1+\varepsilon})\right)={\cal O}(h^\infty):L^2(Y)\to L^2(Y)\eqno{(5.31)}$$
 for every integer $k\ge 0$. Clearly, (5.31) holds trivially for $k=0$. It is easy also to see that (5.31) with $k=1$ implies
 (5.31) for every $k\ge 1$. On the other hand, since 
 $$Z{\rm Op}_h\left(\phi(\eta_1|\mu|/h^{1+\varepsilon})\right)=
 {\rm Op}_h\left((A(0)-a_0)\phi(\eta_1|\mu|/h^{1+\varepsilon})\right)$$
 and $\phi_1=1$ on supp$\,\phi$, (5.31) with $k=1$ follows from Proposition 4.2.
 \eproof
 
 Thus, by Proposition 5.5 we get that the parametrix $\widetilde u_1$ has the following properties.
 
 \begin{Theorem} For all $s\ge 0$, we have the bounds
 $$\left\|P_0\widetilde u_1\right\|_{H^s({\bf R}^+\times Y)}\le C_{s,M}h^{M\varepsilon/2}\|f\|_{L^2(Y)},\eqno{(5.32)}$$
  $$\left\|\widetilde u_1|_{t=0}-{\rm Op}_h\left(\phi(\eta_1|\mu|/h^{1+\varepsilon})\right)f\right\|_{L^2(Y)}\le {\cal O}(h^{\infty})\|f\|_{L^2(Y)},\eqno{(5.33)}$$
   $$\left\|{\cal D}_t\widetilde u_1|_{t=0}\right\|_{L^2(Y)}\le Ch^{\varepsilon}\|f\|_{L^2(Y)}.\eqno{(5.34)}$$
 \end{Theorem}
 
\section{Parametrix construction in the region ${\cal G}_2(\varepsilon)$}

In this section we will construct a parametrix, $\widetilde u_2$, of the solution of the equation
(5.1) with $f_1={\rm Op}_h(\phi_2(\eta_1))f$, where $\phi_2\in C_0^\infty({\bf R})$ is such that on supp$\,\phi_2$ we have
$$|\mu|\sqrt{|\mu|+|\eta_1|}\ge h^{1-\varepsilon},\eqno{(6.1)}$$
$$|\mu|+|\eta_1|\le {\cal O}(h^{\varepsilon}).\eqno{(6.2)}$$
Let $\rho$ be the solution to the equation
$$\rho^2+\eta_1+i\mu q(y,\eta)=0$$
with ${\rm Im}\,\rho>0$. We will be looking for $\widetilde u_2$ in the form
$$\widetilde u_2={\rm Op}_h\left(A(t)\right)f,$$
$$A(t)=\phi(t/|\rho|^2\delta_1)a(t,y,\eta;\mu,h)e^{i\varphi(t,y,\eta;\mu)/h},$$
where $\phi$ is the same function as in the previous section, $\delta_1>0$ is a small constant to be fixed later on, $a=\phi_2(\eta_1)$, $\varphi=0$ for $t=0$. The phase $\varphi$ is independent of $h$ and is of the form
$$\varphi=\sum_{k=1}^{M}t^k\varphi_k$$
where $\varphi_k$ do not depend on $t$, $M\gg 1$ being an arbitrary but fixed integer. The amplitude $a$ is of the form
$$a=\sum_{0\le k+\nu\le M}h^kt^\nu a_{k,\nu}$$
where the functions $a_{k,\nu}$ do not depend on $t$. Note that the identity (5.3) still holds with the new function $A=
\phi(t/|\rho|^2\delta_1)e^{i\varphi/h}a$. Moreover, we have the identity
$$e^{-i\varphi/h}({\cal D}_t^2+t+\eta_1+i\mu q(y,\eta)-ih\partial_{y_1})(e^{i\varphi/h}a)$$ 
$$=-2ih\partial_t\varphi\partial_ta-h^2\partial_t^2a-ih\partial_{y_1}a+((\partial_t\varphi)^2+\partial_{y_1}\varphi+t-\rho^2)a$$
$$=-2ih\sum_{0\le k+\nu\le 2M-2}h^kt^\nu \sum_{j=0}^\nu (j+1)(\nu+1-j)\varphi_{\nu+1-j}\,a_{k,j+1}$$
$$-h\sum_{0\le k+\nu\le M-1}(\nu+1)(\nu+2)h^kt^\nu a_{k-1,\nu+2}-ih\sum_{0\le k+\nu\le M}h^kt^\nu \partial_{y_1}a_{k,\nu}$$ $$+((\partial_t\varphi)^2+\partial_{y_1}\varphi+t-\rho^2)a.\eqno{(6.3)}$$
Let $E_j(t)$, ${\cal E}_j(t)$, $j=1,2$ be defined as in the previous section with the new $A$.
Given a multi-index $\alpha=(\alpha_1,...,\alpha_{n-1})$, set
$$g_\alpha(\varphi)=\lim_{h\to 0}\frac{(-ih)^{|\alpha|}}{|\alpha|!}e^{-i\varphi/h}\partial_y^\alpha(e^{i\varphi/h})=\frac{1}{|\alpha|!}\prod_{j=1}^{n-1}(\partial_{y_j}\varphi)^{\alpha_j}.$$  
The phase satisfies the eikonal equation
$$(\partial_t\varphi)^2+\partial_{y_1}\varphi+t-\rho^2+i\mu\sum_{|\alpha|=1}^M g_\alpha(\varphi)=R_M(t)\eqno{(6.4)}$$
where $R_M(t)={\cal O}(t^{M+1})$ as $t\to 0$. It is easy to see that we have the identities
$$(\partial_t\varphi)^2=\sum_{K=0}^{2M}t^K\sum_{k+j=K}(k+1)(j+1)\varphi_{k+1}\varphi_{j+1},$$
$$\sum_{|\alpha|=1}^M g_\alpha(\varphi)=\sum_{K=1}^{M^2}t^K\sum_{j=1}^M\sum_{k_i\ge 1,k_1+...+k_j=K}\sum_{|\alpha_i|=1}\gamma_{\alpha_1,...,\alpha_j,k_1,...,k_j}\partial_y^{\alpha_1}\varphi_{k_1}...
\partial_y^{\alpha_j}\varphi_{k_j}$$
where $\gamma_{\alpha_1,...,\alpha_j,k_1,...,k_j}$ are constants. 
Thus, if we choose $\varphi_k$ satisfying the equations
$$\varphi_1^2-\rho^2=0, \eqno{(6.5)}$$
$$\sum_{k+j=K}(k+1)(j+1)\varphi_{k+1}\varphi_{j+1}+\partial_{y_1}\varphi_K+\epsilon_K$$ 
$$=-i\mu\sum_{j=1}^M\sum_{k_i\ge 1,k_1+...+k_j=K}\sum_{|\alpha_i|=1}\gamma_{\alpha_1,...,\alpha_j,k_1,...,k_j}\partial_y^{\alpha_1}\varphi_{k_1}...
\partial_y^{\alpha_j}\varphi_{k_j},\quad K\ge 1, \eqno{(6.6)}$$
where $\epsilon_1=1$, $\epsilon_K=0$ for $K\ge 2$, 
then $\varphi$ satisfies the equation (6.4) with
$$R_M(t)=\sum_{K=M+1}^{2M}t^{K}\sum_{k+j=K}(k+1)(j+1)\varphi_{k+1}\varphi_{j+1}$$
$$+i\mu\sum_{K=M+1}^{M^2}t^K\sum_{j=1}^M\sum_{k_i\ge 1,k_1+...+k_j=K}\sum_{|\alpha_i|=1}\gamma_{\alpha_1,...,\alpha_j,k_1,...,k_j}\partial_y^{\alpha_1}\varphi_{k_1}...
\partial_y^{\alpha_j}\varphi_{k_j}.$$
Clearly, $\varphi_1=\rho$ is a solution of (6.5). Then, given $\varphi_1$, ..., $\varphi_K$, $K\ge 1$, we can determine $\varphi_{K+1}$ uniquely
from (6.6).

\begin{lemma} For all integers $k\ge 2$ and all multi-indices $\alpha$ we have the bounds
$$|\partial_y^\alpha\varphi_k|\le C_{k,\alpha}|\rho|^{3-2k},\eqno{(6.7)}$$
$$|{\rm Im}\,\partial_y^\alpha\varphi_k|\le C_{k,\alpha}|\rho|^{2-2k}{\rm Im}\,\rho.\eqno{(6.8)}$$
We also have the bound
$$|\partial_y^\alpha(|\rho|^{-2})|\le C_{\alpha}|\rho|^{-2}.\eqno{(6.9)}$$
Moreover, if $0<t\le \delta_1|\rho|^2$ with a constant $\delta_1>0$ small enough, we have
$${\rm Im}\,\varphi\ge t\,{\rm Im}\,\rho/2 .\eqno{(6.10)}$$
\end{lemma}

{\it Proof.} The bound (6.7) with $k=1$ follows easily by induction in $|\alpha|$ from the identity
$$\sum_{|\alpha_1|+|\alpha_2|=|\alpha|}\gamma_{\alpha_1,\alpha_2}\partial_y^{\alpha_1}\rho\partial_y^{\alpha_2}\rho=i\mu
\partial_y^\alpha q(y,\eta)$$
for $|\alpha|\ge 1$, $\gamma_{\alpha_1,\alpha_2}\neq 0$ being some constants, together with the fact that $\mu={\cal O}(|\rho|^2)$.
The proof of (6.9) is similar, using that
$$|\rho|^2=\eta_1^2+\mu^2q(y,\eta)^2$$
together with the identity
$$\sum_{|\alpha_1|+|\alpha_2|=|\alpha|}\gamma_{\alpha_1,\alpha_2}\partial_y^{\alpha_1}(|\rho|^{-2})\partial_y^{\alpha_2}(|\rho|^2)=0$$
for $|\alpha|\ge 1$. To prove (6.7) for all $k\ge 2$ and all multi-indices $\alpha$ we will proceed by induction in $k+|\alpha|$. Suppose first
that (6.7) holds for all $k\le K$. Then the right-hand side of (6.6) is $\sum_{j=1}^M{\cal O}(|\rho|^{3j-2K})={\cal O}(|\rho|^{3-2K})$.
Thus by (6.6) we get that $\rho\varphi_{K+1}={\cal O}(|\rho|^{2-2K})$, which is the desired bound for $\varphi_{K+1}$. To bound
$\partial_y^\alpha\varphi_{K+1}$ we apply the operator $\partial_y^\alpha$ to the equation (6.6) and proceed in the same way.
The proof of (6.8) is similar, using that $|\mu|\le C|\rho|{\rm Im}\,\rho$ together with the inequality
$$|{\rm Im}\,(z_1...z_k)|\le C_k|z_1|...|z_k|\sum_{j=1}^k\frac{|{\rm Im}\,z_j|}{|z_j|}.$$
To prove (6.10) we use (6.8) to obtain, for $0<t\le \delta_1|\rho|^2$, 
$${\rm Im}\,\varphi=\sum_{k=1}^{M}t^k{\rm Im}\,\varphi_{k}\ge t\,{\rm Im}\,\rho\left(1-C\sum_{k=0}^{M-1}t^k|\rho|^{-2k}\right)$$
$$\ge t\,{\rm Im}\,\rho\left(1-{\cal O}(\delta_1)\right)\ge t\,{\rm Im}\,\rho/2 $$
provided $\delta_1$ is taken small enough. 
\eproof

Set
$$\widetilde E_1(t)=\frac{i\mu}{h}\sum_{|\alpha|=1}^M\frac{(-ih)^{|\alpha|}}{|\alpha|!}
\partial_\eta^\alpha q\left(e^{-i\varphi/h}\partial_y^\alpha(e^{i\varphi/h}a)-g_\alpha(\varphi)a\right),$$
$$\widetilde E_2(t)=\sum_{|\alpha|=0}^M\frac{(-ih)^{|\alpha|}}{|\alpha|!}
\partial_\eta^\alpha \widetilde qe^{-i\varphi/h}\partial_y^\alpha (e^{i\varphi/h}a).$$

\begin{lemma} We have the identities
$$\widetilde E_j(t)=\sum_{k+\nu\le M(M+1)}h^kt^\nu\sum_{|\alpha|=0}^{M}\sum_{k'=0}^k\sum_{\nu'=0}^\nu\widetilde b_{\alpha,k,k',\nu,\nu'}^{(j)}\partial_y^\alpha a_{k',\nu'}\eqno{(6.11)}$$
where the functions $\widetilde b_{\alpha,k,k',\nu,\nu'}^{(j)}$ do not depend on $t$, $h$ and the functions $a_{k,\nu}$, and satisfy the bounds
$$\left|\partial_y^\beta \widetilde b_{\alpha,k,k',\nu,\nu'}^{(j)}\right|\le C_\beta|\rho|^{-2\nu+2\nu'}\eqno{(6.12)}$$
for every multi-index $\beta$. 
\end{lemma}

{\it Proof.} We will first prove by induction in $|\alpha|$ the identity
$$e^{-i\varphi/h}(-ih\partial_y)^\alpha (e^{i\varphi/h})=\sum_{k=0}^{|\alpha|}\sum_{\nu=0}^{M|\alpha|}h^kt^\nu c_{\alpha,k,\nu}\eqno{(6.13)}$$
with functions $c_{\alpha,k,\nu}$ independent of $t$, $h$ and satisfying the bounds
$$\left|\partial_y^\beta c_{\alpha,k,\nu}\right|\le C_\beta|\rho|^{-2\nu}\eqno{(6.14)}$$
for every multi-index $\beta$. Let $\alpha=\alpha_1+\alpha_2$ with $|\alpha_1|=1$ and suppose (6.13) fulfilled with $\alpha_2$. Then we have
$$e^{-i\varphi/h}(-ih\partial_y)^{\alpha_1+\alpha_2}(e^{i\varphi/h})=e^{-i\varphi/h}(-ih\partial_y)^{\alpha_1}e^{i\varphi/h}
\sum_{k=0}^{|\alpha_2|}\sum_{\nu=0}^{M|\alpha_2|}h^kt^\nu c_{\alpha_2,k,\nu}$$
$$=\partial_y^{\alpha_1}\varphi\sum_{k=0}^{|\alpha_2|}\sum_{\nu=0}^{M|\alpha_2|}h^kt^\nu c_{\alpha_2,k,\nu}-i
\sum_{k=0}^{|\alpha_2|}\sum_{\nu=0}^{M|\alpha_2|}h^{k+1}t^\nu \partial_y^{\alpha_1}c_{\alpha_2,k,\nu}$$
$$=\sum_{k=0}^{|\alpha_2|}\sum_{\nu=0}^{M|\alpha_2|+M}h^kt^\nu \sum_{\ell=1}^\nu\partial_y^{\alpha_1}\varphi_\ell\, c_{\alpha_2,k,\nu-\ell}-i
\sum_{k=0}^{|\alpha_2|+1}\sum_{\nu=0}^{M|\alpha_2|}h^{k}t^\nu \partial_y^{\alpha_1}c_{\alpha_2,k-1,\nu}.$$
Hence (6.13) holds for $\alpha_1+\alpha_2$ with
$$c_{\alpha_1+\alpha_2,k,\nu}=\sum_{\ell=1}^\nu\partial_y^{\alpha_1}\varphi_\ell\, c_{\alpha_2,k,\nu-\ell}-i\partial_y^{\alpha_1}c_{\alpha_2,k-1,\nu}.\eqno{(6.15)}$$
It follows from (6.7) and (6.15) that if (6.14) holds with $\alpha_2$, it holds  with $\alpha_1+\alpha_2$, which proves the assertion.

Using (6.13) we can write
$$e^{-i\varphi/h}(-ih\partial_y)^\alpha (e^{i\varphi/h}a)=\sum_{|\alpha_1|+|\alpha_2|=|\alpha|}\gamma_{\alpha_1,\alpha_2}
e^{-i\varphi/h}(-ih\partial_y)^{\alpha_1} (e^{i\varphi/h})(-ih\partial_y)^{\alpha_2}a$$ 
$$=\sum_{k=0}^{|\alpha|}\sum_{\nu=0}^{M|\alpha|}h^kt^\nu \sum_{|\alpha_1|+|\alpha_2|=|\alpha|}\gamma_{\alpha_1,\alpha_2}c_{\alpha_1,k-|\alpha_2|,\nu}(-i\partial_y)^{\alpha_2}a.$$
It follows from this identity and (6.14) that the functions $\widetilde E_j$ are of the form
$$\widetilde E_j(t)=\sum_{k=0}^{M}\sum_{\nu=0}^{M^2}\sum_{|\alpha|=0}^M h^kt^\nu\widetilde c_{\alpha,k,\nu}^{(j)}\partial_y^{\alpha}a
\eqno{(6.16)}$$
with functions $\widetilde c_{\alpha,k,\nu}^{(j)}$ independent of $t$, $h$ and $a$, and satisfying the bounds 
$\partial_y^\beta \widetilde c_{\alpha,k,\nu}^{(j)}={\cal O}_\beta|(\rho|^{-2\nu})$, $\forall\beta$.
Now (6.11) follows from (6.16) with
$$\widetilde b_{\alpha,k,k',\nu,\nu'}^{(j)}=\widetilde c_{\alpha,k-k',\nu-\nu'}^{(j)}.$$
\eproof

We let now the functions $a_{k,\nu}$ satisfy the equations
$$2i\sum_{j=0}^\nu (j+1)(\nu+1-j)\varphi_{\nu+1-j}\,a_{k,j+1}+(\nu+1)(\nu+2)a_{k-1,\nu+2}+i\partial_{y_1}a_{k,\nu}$$
$$=\sum_{j=1}^2\sum_{|\alpha|=0}^{M}\sum_{k'=0}^k\sum_{\nu'=0}^\nu\widetilde b_{\alpha,k,k',\nu,\nu'}^{(j)}\partial_y^\alpha a_{k',\nu'},\eqno{(6.17)}$$
$a_{0,0}=\phi_2(\eta_1)$, $a_{k,0}=0$ for $k\ge 1$, $a_{-1,\nu}=0$, $\nu\ge 0$. Let $K,J\ge 0$ be any integers. Now it is clear that, given
$a_{k,\nu}$ for $k\le K$, $\forall \nu\ge 0$, and $a_{K+1,\nu}$ for $\nu\le J$, we can determine $a_{K+1,J+1}$ from (6.17). Therefore,
by (6.17) we can find all $a_{k,\nu}$. Moreover, using (6.7) and (6.12) one can easily prove the following

\begin{lemma} For all integers $k,\nu\ge 0$ and all multi-indices $\alpha$ we have the bounds
$$\left|\partial_y^\alpha a_{k,\nu}\right|\le C_{k,\nu,\alpha}|\rho|^{-3k-2\nu}.\eqno{(6.18)}$$
\end{lemma}

In view of (6.3) and (6.11), in this case we still have the identity (5.17) with a function $B$ of the form
$$B(t)=e^{i\varphi/h}\phi(t/|\rho|^2\delta_1)B_1(t)+B_2(t),$$
where 
$$B_1(t)=-2ih\sum_{M+1\le k+\nu\le 2M-2}h^kt^\nu \sum_{j=0}^\nu (j+1)(\nu+1-j)\varphi_{\nu+1-j}\,a_{k,j+1}$$
$$+h\sum_{k+\nu = M}(\nu+1)(\nu+2)h^kt^\nu a_{k-1,\nu+2}+R_M(t)a$$
$$+\sum_{j=1}^2\sum_{M+1\le k+\nu\le M(M+1)}h^kt^\nu\sum_{|\alpha|=0}^{M}\sum_{k'=0}^k\sum_{\nu'=0}^\nu\widetilde b_{\alpha,k,k',\nu,\nu'}^{(j)}\partial_y^\alpha a_{k',\nu'},$$
$$B_2(t)=\left[{\cal D}_t^2-ih\partial_{y_1},\phi(t/|\rho|^2\delta_1)\right]e^{i\varphi/h}a$$
$$+\frac{i\mu}{h}\sum_{|\alpha|=1}^M\frac{(-ih)^{|\alpha|}}{|\alpha|!}\partial_\eta^\alpha q\left(\partial_y^\alpha(\phi e^{i\varphi/h}a)-
\phi\partial_y^\alpha(e^{i\varphi/h}a)\right)$$
$$+\sum_{|\alpha|=0}^M\frac{(-ih)^{|\alpha|}}{|\alpha|!}\partial_\eta^\alpha \widetilde q\left(\partial_y^\alpha(\phi e^{i\varphi/h}a)-
\phi\partial_y^\alpha(e^{i\varphi/h}a)\right).$$
 Lemmas 6.1 and 6.3 imply the following 
 
 \begin{lemma} For all multi-indices $\alpha$ we have the bounds
 $$\left|\partial_y^\alpha B(t)\right|\le C_{\alpha}h^{\varepsilon M-|\alpha|},\eqno{(6.19)}$$
 $$\left|\partial_y^\alpha A(t)\right|\le C_{\alpha}h^{-(1-3\varepsilon )|\alpha|}.\eqno{(6.20)}$$
 \end{lemma}
 
 {\it Proof}. Note first that the condition (6.1) implies
 $$\frac{h}{|\rho|^3}\le\frac{ C_1 h}{|\mu||\rho|}\le C_2h^\varepsilon\eqno{(6.21)}$$
 with some constants $C_1,C_2>0$. By (6.7), (6.18) and (6.21) we have, for $0\le t\le\delta_1|\rho|^2$,
 $$h^kt^\nu\left|e^{i\varphi/h}a_{k,\nu}\right|\le C_{k,\nu}\left(\frac{h}{|\rho|^3}\right)^k\left(\frac{t}{|\rho|^2}\right)^\nu e^{-t{\rm Im}\,\rho/2h}$$
 $$\le C_{k,\nu}\left(\frac{h}{|\rho|^3}\right)^k\left(\frac{h}{|\mu||\rho|}\right)^\nu\le C_{k,\nu}h^{\varepsilon(k+\nu)}\eqno{(6.22)}$$
 where we have used that $|\rho|{\rm Im}\,\rho\ge C|\mu|$ with some constant $C>0$. In the same way, since
 $e^{-i\varphi/h}(h\partial_y)^\alpha(e^{i\varphi/h})={\cal O}_\alpha(1)$ for $0\le t\le 1$, one can get
 that for any multi-index $\alpha$ and for $0\le t\le\delta_1|\rho|^2$,
  $$h^kt^\nu\left|(h\partial_y)^\alpha\left(e^{i\varphi/h}a_{k,\nu}\right)\right|\le C_{\alpha,k,\nu}h^{\varepsilon(k+\nu)}.\eqno{(6.23)}$$
  It follows easily from (6.23) that, for $0\le t\le\delta_1|\rho|^2$,
  $$\left|(h\partial_y)^\alpha\left(e^{i\varphi/h}B_1(t)\right)\right|\le C_{\alpha}h^{\varepsilon M}.\eqno{(6.24)}$$
 On the other hand, for $\frac{\delta_1}{2}|\rho|^2\le t\le\delta_1|\rho|^2$, we have
 $$\left|e^{i\varphi/h}\right|\le e^{-\delta_1|\rho|^2{\rm Im}\,\rho/4h}\le e^{-c_1|\rho||\mu|/h}\le e^{-c_2h^{-\varepsilon}}
 \eqno{(6.25)}$$
 with some constants $c_1,c_2>0$. In view of (6.9) we have $\partial_y^\alpha\phi(t/|\rho|^2\delta_1)={\cal O}_\alpha(1)$, $\forall\alpha$,
 and $\partial_t^\ell\phi(t/|\rho|^2\delta_1)={\cal O}_\ell(|\mu|^{-\ell})={\cal O}_\ell(h^{-\ell})$, $\forall\ell$.
 Therefore, by (6.23) and (6.25) we obtain
  $$\left|\partial_y^\alpha B_2(t)\right|\le C_{\alpha}e^{-ch^{-\varepsilon}}\eqno{(6.26)}$$
 with some constant $c>0$. Thus (6.19) follows from (6.24) and (6.26).
 
 To prove (6.20) we need to improve the estimate (6.23) when $|\alpha|\ge 1$. To this end, observe that by Lemma 6.1
 we have $\partial_y^\alpha\varphi={\cal O}_\alpha(t|\rho|)={\cal O}_\alpha(|\rho|^3)$, $\forall\alpha$, for $0\le t\le\delta_1|\rho|^2$.
 Therefore, by induction in $|\alpha|$ one easily gets
 $$\left|e^{-i\varphi/h}\partial_y^\alpha(e^{i\varphi/h})\right|\le C_\alpha\left(\frac{|\rho|^3}{h}\right)^{|\alpha|}+C_\alpha.\eqno{(6.27)}$$
 By (6.2), (6.10) and (6.27), for $0\le t\le\delta_1|\rho|^2$,
 $$\left|\partial_y^\alpha(e^{i\varphi/h})\right|\le C_\alpha\left(\frac{|\rho|^3}{h}\right)^{|\alpha|}+C_\alpha \le C_\alpha h^{-(1-3\varepsilon)|\alpha|}.\eqno{(6.28)}$$
 On the other hand, by (6.18) we have $\partial_y^\alpha a={\cal O}_\alpha(1)$ for $0\le t\le\delta_1|\rho|^2$. Therefore, (6.20) follows from
 (6.28). 
 \eproof
  
  Lemma 6.4 implies the following

\begin{prop} For all $s\ge 0$, we have the bounds
 $$\left\|P_0\widetilde u_2\right\|_{H^s({\bf R}^+\times Y)}\le C_{s,M}h^{M\varepsilon/2}\|f\|_{L^2(Y)},\eqno{(6.29)}$$
   $$\left\|{\cal D}_t\widetilde u_2|_{t=0}\right\|_{L^2(Y)}\le Ch^{\varepsilon}\|f\|_{L^2(Y)}.\eqno{(6.30)}$$
 \end{prop}

{\it Proof.} By Proposition 4.1 and (6.19), there is $\ell>0$ dpending only on the dimension such that ${\rm Op}_h({\cal D}_y^\alpha 
{\cal D}_t^\beta B(t))={\cal O}_{\alpha,\beta}(h^{M\varepsilon-\ell}):L^2(Y)\to L^2(Y)$, $\forall\alpha,\beta$, while
by Proposition 4.2 and (6.20) we have ${\cal D}_y^\alpha 
{\cal D}_t^\beta {\cal E}_j(t)={\cal O}_{\alpha,\beta}(h^{M\varepsilon-\ell}):L^2(Y)\to L^2(Y)$, $\forall\alpha,\beta$. 
This implies (6.29) in view of the identity (5.17).

To prove (6.30), observe that
$${\cal D}_t\widetilde u_2|_{t=0}={\rm Op}_h\left(\rho-ih\sum_{k=0}^{M-1}h^ka_{k,1}\right)f.$$
In view of (6.2) and (6.7), we have $\partial_y^\alpha\rho={\cal O}_\alpha(|\rho|)={\cal O}_\alpha(h^\varepsilon)$, and hence
by Proposition 4.1 we get ${\rm Op}_h(\rho)={\cal O}_\alpha(h^\varepsilon):L^2(Y)\to L^2(Y)$. Furthermore, by
(6.18) we also have $h^{k+1}\partial_y^\alpha a_{k,1}={\cal O}_\alpha(|\rho|)={\cal O}_\alpha(h^\varepsilon)$, and we apply
once again Proposition 4.1 to get (6.30).
\eproof

To complete the construction of our parametrix $\widetilde u$ we will consider two cases.

Case 1. $h^{(1+\varepsilon)/2}\le|\mu|\le h^{\varepsilon}$, $0<\varepsilon\ll 1$. Then the condition (6.1) is fulfilled for all
$\eta_1$. We take $\widetilde u=\widetilde u_2$, where $\widetilde u_2$ is the parametrix constructed above with $\phi_2(\eta_1)=
\phi(\eta_1/h^\varepsilon)$. Clearly the condition (6.2) is fullfiled as long as $\eta_1\in{\rm supp}\,\phi_2$.

Case 2. $h^{1-2\varepsilon}\le|\mu|\le h^{(1+\varepsilon)/2}$. Then $(\mu,\eta_1)\in{\cal G}_1(\varepsilon)$ as long as 
$\eta_1\in{\rm supp}\,\phi(\eta_1|\mu|/h^{1+\varepsilon})$. We take $\widetilde u=\widetilde u_1+\widetilde u_2$, where 
$\widetilde u_1$ is the parametrix constructed in Section 5 and 
$\widetilde u_2$ is the parametrix constructed in Section 6 with $\phi_2(\eta_1)=\phi(\eta_1/h^\varepsilon)-
\phi(\eta_1|\mu|/h^{1+\varepsilon})$. Clearly $\eta_1={\cal O}(h^\varepsilon)$ on supp$\,\phi_2$, and hence the condition (6.2) is
fulfilled in this case. Moreover, if $(\mu,\eta_1)\in {\cal G}_2(\varepsilon)$, then
$$|\mu|\sqrt{|\mu|+|\eta_1|}\ge |\mu|^{1/2}h^{(1+\varepsilon)/2}\ge h^{1-\varepsilon/2}.$$
Hence, with this choice of the function $\phi_2$, the condition (6.1) is satisfied (with $\varepsilon/2$
in place of $\varepsilon$) as long as $\eta_1\in{\rm supp}\,\phi_2$. 

In both cases the operator $\widetilde N$ defined by $\widetilde Nf:={\cal D}_t\widetilde u|_{t=0}$
provides a parametrix for the DN map $f\to {\cal D}_t u|_{t=0}$, where $u$ is the solution to the equation (5.1)
with $u|_{t=0}={\rm Op}_h(\phi(\eta_1/h^\varepsilon))f$. It follows from Theorem 5.7 and Proposition 6.5 that
$\widetilde u$ and $\widetilde N$ have the following properties.

 \begin{Theorem} For all $s\ge 0$, we have the bounds
 $$\left\|P_0\widetilde u\right\|_{H^s({\bf R}^+\times Y)}\le C_{s,M}h^{M\varepsilon/2}\|f\|_{L^2(Y)},\eqno{(6.31)}$$
  $$\left\|\widetilde u|_{t=0}-{\rm Op}_h\left(\phi(\eta_1/h^{\varepsilon})\right)f\right\|_{L^2(Y)}\le {\cal O}(h^{\infty})\|f\|_{L^2(Y)},\eqno{(6.32)}$$
   $$\left\|\widetilde Nf\right\|_{L^2(Y)}\le Ch^{\varepsilon/2}\|f\|_{L^2(Y)},\eqno{(6.33)}$$
   $$\left\|{\rm Op}_h((1-\phi_1)(\eta_1/h^\varepsilon))\widetilde Nf\right\|_{L^2(Y)}\le {\cal O}(h^{\infty})\|f\|_{L^2(Y)},\eqno{(6.34)}$$
   where $\phi_1\in C_0^\infty({\bf R})$ is independent of $h$ and $\mu$, and $\phi_1=1$ on ${\rm supp}\,\phi$.
 \end{Theorem}
 
 Note that the estimate (6.34) follows from Proposition 4.2 in the same way as in the proof of Lemma 5.6.

\section{Eigenvalue-free regions}

In this section we will study the problem
$$\left\{
\begin{array}{lll}
\left(h^2\nabla c_1(x)\nabla + zn_1(x)\right)u_1=0 &\mbox{in} &\Omega,\\
\left(h^2\nabla c_2(x)\nabla + zn_2(x)\right)u_2=0 &\mbox{in} &\Omega,\\
u_1=u_2,\,\,\, c_1\partial_\nu u_1=c_2\partial_\nu u_2& \mbox{on}& \Gamma,
\end{array}
\right.
\eqno{(7.1)}
$$
where $0<h\ll 1$, $z=1+i\,{\rm Im}\,z$, $0<|{\rm Im}\,z|\le 1$. Denote by $N_j(h,z)$, $j=1,2$, the Dirichlet-to-Neumann map
corresponding to the Laplacian $n_j(x)^{-1}\nabla c_j(x)\nabla$ introduced in Section 2 (with $\mu={\rm Im}\,z$). In this section
we will prove the following

\begin{Theorem} Under the conditions of Theorem 1.1, given any $0<\varepsilon\ll 1$ there is $h_0(\varepsilon)>0$ so that
the operator
$$T(h,z)=c_1N_1(h,z)-c_2N_2(h,z):H^1(\Gamma)\to L^2(\Gamma)$$
is invertible for $0<h\le h_0$, $|{\rm Im}\,z|\ge h^{1-\varepsilon}$.
\end{Theorem}

{\it Proof.} We may suppose that $|{\rm Im}\,z|\le h^{\varepsilon}$ since for $h^{\varepsilon}\le |{\rm Im}\,z|\le 1$ the theorem is proved in
\cite{kn:V}. 
Let $\Delta_\Gamma$ be the negative Laplace-Beltrami operator on $\Gamma$ with the Riemannian metric induced by the Euclidean one in ${\bf R}^d$. Denote by $r_0(x',\xi')$ the principal symbol of $-\Delta_\Gamma$ written in the coordinates $(x',\xi')\in T^*\Gamma$.
Set $\Sigma_j(\varepsilon)=\left\{(x',\xi')\in T^*\Gamma:|r_0-m_j|\le h^{\varepsilon/2}\right\}$, 
where $m_j$ denotes the restriction on $\Gamma$ of the function $n_j/c_j$. It is easy to see that
the conditions (1.2) and (1.3) imply $\Sigma_1(\varepsilon)\cap\Sigma_2(\varepsilon)=\emptyset$, provided $h$ is taken small enough.
Throughout this section, $\rho_j$, $j=1,2$, will denote the solution to the equation
$$\rho^2+r_0(x',\xi')-zm_j(x')=0$$
 with ${\rm Im}\,\rho>0$. Observe  that
$$c_1\rho_1-c_2\rho_2=\frac{\widetilde c(x')(c_0(x')r_0(x',\xi')-z)}{c_1\rho_1+c_2\rho_2}\eqno{(7.2)}$$
where $\widetilde c$ and $c_0$ are the restrictions on $\Gamma$ of the functions
$$c_1n_1-c_2n_2\quad\mbox{and}\quad\frac{c_1^2- c_2^2}{c_1n_1-c_2n_2}$$
respectively. Clearly, under the conditions of Theorem 1.1, we have $\widetilde c(x')\neq 0$, $\forall x'\in \Gamma$.
Moreover, (1.2) implies $c_0\equiv 0$, while (1.3) implies $c_0(x')<0$, $\forall x'\in \Gamma$.
 Hence, under the conditions of Theorem 1.1, we have $c_1^2\rho_1^2\neq c_2^2\rho_2^2$ on $\Gamma$ as $|{\rm Im}\,z|\to 0$.
 It is easy to see that $|\rho_j|\ge Const>0$ on $\Sigma_{3-j}(\varepsilon)$, $j=1,2$. Let $\chi_\varepsilon^{(j)}\in C_0^\infty(T^*\Gamma)\cap
 {\cal S}^0_{\varepsilon/2}$, $\chi_\varepsilon^{(j)}=1$ on $\Sigma_j(\varepsilon)$, $\chi_\varepsilon^{(j)}=0$ outside a larger
 ${\cal O}(h^{\varepsilon/2})$ neighbourhood of $\{r_0=m_j\}$. Then we have $\widetilde\rho_j=(1-\chi_\varepsilon^{(j)})\rho_j\in {\cal S}^1_{\varepsilon/2}$.
 
 By (7.2) we also have
 $$C_1\langle r_0\rangle^{k/2}\le |c_1\rho_1-c_2\rho_2|\le C_2\langle r_0\rangle^{k/2},\quad C_2>C_1>0,\eqno{(7.3)}$$
  where $k=-1$ if (1.2) holds, $k=1$ if (1.3) holds. Since $\chi_\varepsilon^{(j)}\rho_j={\cal O}(h^{\varepsilon/4})$, (7.3) remains valid with $\widetilde\rho_j$ in place of
  $\rho_j$. Using this we will prove the following
 
 \begin{prop} Under the conditions of Theorem 1.1, we have the estimate
 $$\left\|T(h,z)f-{\rm Op}_h(c_1\widetilde\rho_1-c_2\widetilde\rho_2)f\right\|_{H^{\frac{1-k}{2}}(\Gamma)}\le Ch^{\varepsilon/4}\|f\|_{H^{\frac{k+1}{2}}(\Gamma)}\eqno{(7.4)}$$
 for $0<h\le h_0$, $|{\rm Im}\,z|\ge h^{1-\varepsilon}$.
 \end{prop}
 
 {\it Proof.} Let $\chi\in C_0^\infty(T^*\Gamma)$, $\chi=1$ on $\{r_0\le R_0\}$ with some constant $R_0\gg 1$. The estimate (7.4)
 with $f$ replaced by ${\rm Op}_h(1-\chi)f$ is proved, under the conditions (1.2) and (1.3), in Section 5 of \cite{kn:V} (see also \cite{kn:PV}).
 Therefore, to prove (7.4) it suffices to show that
 $$\left\|N_j(h,z){\rm Op}_h(\chi)f-{\rm Op}_h(c_j\widetilde\rho_j\chi)f\right\|_{L^2(\Gamma)}\le Ch^{\varepsilon/4}\|f\|_{L^2(\Gamma)}.\eqno{(7.5)}$$
Let $\widetilde\chi_\varepsilon^{(j)}\in C_0^\infty(T^*\Gamma)\cap
 {\cal S}^0_{\varepsilon/2}$, $\widetilde \chi_\varepsilon^{(j)}=1$ on $\Sigma_j(\varepsilon)$, be such that  $\widetilde \chi_\varepsilon^{(j)}\widetilde\rho_j\equiv 0$. Then (7.5) with $\widetilde \chi_\varepsilon^{(j)}$ in place of $\chi$
 follows from the estimate (2.7) of Theorem 2.2, while (7.5) with $\chi-\widetilde \chi_\varepsilon^{(1)}-\widetilde \chi_\varepsilon^{(2)}$ in place of $\chi$
 follows from the estimates (2.4) and (2.5) of Theorem 2.1. 
\eproof

Thus we have reduced the problem to that one of inverting the operator ${\cal A}={\rm Op}_h(c_1\widetilde\rho_1-c_2\widetilde\rho_2)$.
This, however, is much easier since the symbol $c_1\widetilde\rho_1-c_2\widetilde\rho_2\in {\cal S}^k_{\varepsilon/2}$ is elliptic in view of 
(7.3). Hence $(c_1\widetilde\rho_1-c_2\widetilde\rho_2)^{-1}\in {\cal S}^{-k}_{\varepsilon/2}$ and there exists an inverse 
${\cal A}^{-1}={\cal O}(1):H^{\frac{1-k}{2}}(\Gamma)\to H^{\frac{1+k}{2}}(\Gamma)$. 
Then (7.4) yields
$$\|f\|_{H^{\frac{k+1}{2}}(\Gamma)}\le C\|{\cal A}^{-1}T(h,z)f\|_{H^{\frac{k+1}{2}}(\Gamma)}+Ch^{\varepsilon/4}\|f\|_{H^{\frac{k+1}{2}}(\Gamma)}$$
which after taking $h$ small enough becomes
$$\|f\|_{H^{\frac{k+1}{2}}(\Gamma)}\le 2C\|{\cal A}^{-1}T(h,z)f\|_{H^{\frac{k+1}{2}}(\Gamma)}.\eqno{(7.6)}$$
Clearly, (7.6) implies the invertibility of the operator $T$ in the desired region.
\eproof

{\bf Acknowledgements.} I would like to thank Vesselin Petkov for some very usefull discussions and suggestions.

G. Vodev

Universit\'e de Nantes,

 D\'epartement de Math\'ematiques, UMR 6629 du CNRS,
 
 2, rue de la Houssini\`ere, BP 92208, 
 
 44332 Nantes Cedex 03, France,
 
 e-mail: vodev@math.univ-nantes.fr

\end{document}